\documentclass[11pt]{article}
\usepackage{amsmath}
\usepackage{amsthm,amsfonts,amsmath,amssymb}
\usepackage{multirow}

\newtheorem{thm}{Theorem}[section]

\newtheorem{lem}[thm]{Lemma}

\numberwithin{equation}{section}
\numberwithin{table}{section}

\def\di{\bigm|} \def\lg{\langle} \def\rg{\rangle}
\def\nd{\mathrel{\bigm|\kern-.7em/}}

\def\f{\noindent}

\def\det{\hbox{\rm det}}

\def\mod{\hbox{\rm mod }}

\def\rank{\hbox{\rm rank}}

\def\demo{\f{\bf Proof}\hskip10pt}

\def\diag{\hbox{\rm diag}}

\def\qed{\hfill $\Box$}

\def\lg{\langle}
\def\rg{\rangle}

\def\rr#1{{\rm (#1)}}

\def\A{\mathcal{A}$$}

\begin{document}
\title{Finite $p$-groups with a minimal non-abelian
subgroup of index $p$ (III)
\thanks{This work was supported by NSFC (no. 11371232),
by NSF of Shanxi Province (no. 2012011001-3 and 2013011001-1). }}
\author{Haipeng Qu, Mingyao Xu and Lijian An
\thanks{Corresponding author. e-mail: anlj@sxnu.edu.cn}\\
Department of Mathematics,
Shanxi Normal University\\
Linfen, Shanxi 041004, P. R. China }

\maketitle

\begin{abstract}
In this paper, we finished the classification of three-generator
finite $p$-groups $G$ such that $\Phi(G)\le Z(G)$. This paper is a
part of classification of finite $p$-groups with a minimal
non-abelian subgroup of index $p$, and partly solved a problem
proposed by Y. Berkovich.

\medskip


\end{abstract}

\baselineskip=16pt

\section{Introduction}
Groups in this paper are all finite $p$-groups. Notation and
terminology are consistent with that in \cite{QYXA} and \cite{ALQZ}.

Let $G$ be a non-abelian $p$-groups.
If $d(G)=2$, then we can prove that $G\in \mathcal{A}_1$ if and only
if $G$ satisfy one of the following conditions:

{\rm (1)} $|G'|=p$;

{\rm (2)} $c(G)=2$ and $G$ has an abelian maximal subgroup;

{\rm (3)} $\Phi(G)\le Z(G)$.

In general, above conditions are not mutually equivalent.
If we replace $d(G)=2$ with $d(G)=3$, then we have
(1)$\Rightarrow$ (2)$ \Rightarrow$(3). Moreover, it is easy to prove
that $G'\le C_p^2$ for {\rm (2)} and $G'\le C_p^3$ for {\rm (3)}.

In \cite{ALQZ}, we completely classified three-generator groups
satisfying condition (2). In this paper, we completely classified
three-generator groups satisfying condition (3) but not satisfying
condition (2). That is, we classified groups $G$ with $d(G)=3$,
$\Phi(G)\le Z(G)$ and $G'\cong C_p^3$. Together with \cite{ALQZ}, we
completely classified groups $G$ with $d(G)=3$ and $\Phi(G)\le
Z(G)$. As a direct application of this classification, we can give a classification about
three-generator groups with at least two $\mathcal{A}_1$-subgroups
of index $p$.

In \cite{Yak}, Y. Berkovich proposed the following

\medskip

{\bf Problem 239:} Classify the $p$-groups containing an
$\mathcal{A}_1$-subgroup of index $p$.

\medskip

In \cite{QYXA}, this problem was divided into two parts:

\medskip

Part 1. Classify the finite $p$-groups with at least two
$\mathcal{A}_1$-subgroups of index $p$.

Part 2. Classify the finite $p$-groups with a unique
$\mathcal{A}_1$-subgroup  of index $p$.

\medskip

Together with \cite{AHZ}, we solve Part 1 of problem 239. Finally,
together with \cite{QYXA, QZGA}, we completely solve problem 239.

By the way, we also look into some properties for those groups we
obtained. In particular, we give the minimal and the maximal index
of $\mathcal{A}_1$-subgroups, and pick out all metahamiltonian groups
(that is, groups whose non-normal subgroups are abelian). These
properties are useful in the classification of
$\mathcal{A}_3$-groups \cite{zhang} and metahamiltonian $p$-groups
\cite{AF}, respectively.

\section{Preliminaries}


%
%
%
In this paper, $p$ is always a prime. We use $F_p$ to denote the
finite field containing $p$ elements. $F_p^*$ is the multiplicative
group of $F_p$. $(F_p^*)^2=\{ a^2\di a\in F_p^*\}$ is a subgroup of
$F_p^*$. $V_m(G)$ means the set of all $p^m$th powers in $G$ and $\mho_m(G)=\lg V_m(G)\rg$. For a finite non-abelian $p$-group $G$, we use $p^{I_{\min}}$ and $p^{I_{\max}}$ to denote the minimal index and the maximal index of $\mathcal{A}_1$-subgroups of $G$ respectively.

We use $\rank(A)$ to denote the rank of $A$, where $A$ is a matrix over $F$. Recall that two matrices $A$ and $B$ over $F$ are called {\it
sub-congruent} if there exists an element $\lambda\in F^*$ and an
invertible matrix $P$ over $F$ such that $\lambda P^tAP=B$.

Suppose that $G$ is a finite $p$-group with $d(G)=3$, $\Phi(G)\le
Z(G)$ and $G'\cong C_p^3$. Let the type of $G/G'$ be
$(p^{m_1},p^{m_2},p^{m_3})$, where $m_1\ge m_2\ge m_3$, and
$$G/G'=\lg a_1G'\rg\times\lg a_2G'\rg\times\lg a_3G'\rg,$$ where
$o(a_iG')=p^{m_i}$, $i=1,2,3$. Then $G=\lg a_1,a_2,a_3\rg$ and
$G'=\lg [a_2,a_3],[a_3,a_1],[a_1,a_2]\rg$. Let $x=[a_2,a_3],
y=[a_3,a_1]$ and $z=[a_1,a_2]$. Since $a_i^{p^{m_i}}\in G'$, we may
assume that $a_i^{p^{m_i}}=x^{w_{i1}}y^{w_{i2}}z^{w_{i3}}$ where
$i=1,2,3$. Then we get a $3\times 3$ matrix $w(G)=(w_{ij})$ over
$F_p$. We call $w(G)$ a characteristic matrix of
$G$. Notice that $w(G)$ will be changed if we change the generators $a_1,a_2,a_3$.  
We also call $a_1,a_2,a_3$ a set of characteristic generators about $w(G)$. In general, $a_1,a_2,a_3$ is always a set of characteristic generators about $w(G)$.

\begin{lem}\label{hetong-p} Suppose that $p$ is odd, $\{1,\eta\}$ is a transversal for
$(F_p^*)^2$ in $F_p^*$. Then the following matrices form a
transversal for the congruence $($sub-congruence$)$ classes of
invertible matrices of order $2$ over $F_p$:

\medskip
\ \ \ {\rm (1)} $\left(
 \begin{array}{cc}
 0& 1\\
 -1& 0
 \end{array}
 \right)$, \ \ \
{\rm (2)} $\left(
 \begin{array}{cc}
 \nu_1  & 1\\
 -1& 0
 \end{array}
 \right)$,\ \ \
{\rm (3)} $\left(
 \begin{array}{cc}
 1& 0\\
 0& \nu_2
 \end{array}
 \right)$,\ \ \
{\rm (4)} $\left(
 \begin{array}{cc}
 1& 1\\
 -1& r
 \end{array}
 \right)$,

\medskip

\noindent where $\nu_1=1$ or $\eta$ $(\nu_1=1)$, $\nu_2=1$ or
$\eta$, $r=1,2,\dots,p-2$.
 \end{lem}

\demo
By \cite[Lemma 4.2]{ALQZ}, matrices (1), (2) and $M=\left(
 \begin{array}{cc}
 1& t\\
 -t& \nu_2
 \end{array}
 \right)$ form a
transversal for the congruence $($sub-congruence$)$ classes of
invertible matrices of order $2$ over $F_p$, where $\nu_2=1$ or
$\eta$, $t\in\{0,1,\dots,\frac{p-1}{2}\}$ such that $t^2\neq
-\nu_2$. If $t=0$, then we get the matrix of Type (3). If $t\neq 0$, then, letting $P=\diag(1,t^{-1})$, we have $P^tMP=\left(
 \begin{array}{cc}
 1& 1\\
 -1& t^{-2}\nu_2
 \end{array}
 \right)$. Hence we get the matrix of Type (4).
\qed

\begin{lem}{\rm (\cite[Lemma 4.3]{ALQZ})}\label{hetong-2}
The following matrices form a transversal for the congruence class of
invertible matrices of order $2$ over $F_2$.

\medskip
 \ \ \ {\rm (1)} $\left(
 \begin{array}{cc}
 1& 0\\
 0& 1
 \end{array}
 \right)$, \ \ \
{\rm (2)} $\left(
 \begin{array}{cc}
 0& 1\\
 1& 0
 \end{array}
 \right)$, \ \ \
{\rm (3)} $\left(
 \begin{array}{cc}
 1& 0\\
 1& 1
 \end{array}
 \right)$.
\end{lem}

\begin{lem}\label{hetong-p-non-invertible} Suppose that $p$ is prime $(p=2$ is possible$)$. For odd $p$, $\{1,\eta\}$ is a transversal for
$(F_p^*)^2$ in $F_p^*$. Then the following matrices form a
transversal for the congruence $($sub-congruence$)$ classes of
non-invertible matrices of order $2$ over $F_p$:

\medskip
\ \ \ {\rm (1)} $\left(
 \begin{array}{cc}
 0& 1\\
 0& 0
 \end{array}
 \right)$, \ \ \
{\rm (2)} $\left(
 \begin{array}{cc}
 0& 0\\
 0& \nu
 \end{array}
 \right)$,\ \ \
{\rm (3)} $\left(
 \begin{array}{cc}
 0& 0\\
 0& 0
 \end{array}
 \right)$,
where $\nu=1$ or $\eta$ $(\nu=1)$.
 \end{lem}

\demo Suppose that $A\ne 0$ be a non-invertible matrix of order $2$
over $F_p$. Then $\rank(A)=1$. Let $A=\left(
 \begin{array}{cc}
 a_{11}& a_{12}\\
 a_{21}& a_{22}
 \end{array}
 \right)$ and $P_1=\left(
 \begin{array}{cc}
 0& 1\\
 1& 0
 \end{array}
 \right)$ . Then $P_1^tAP_1=\left(
 \begin{array}{cc}
 a_{22}& a_{21}\\
 a_{12}& a_{11}
 \end{array}
 \right)$. So, without loss of generality, we may assume that $\left(
 \begin{array}{c}
 a_{12}\\
 a_{22}
 \end{array}
 \right)\ne \left(
 \begin{array}{c}
 0\\
 0
 \end{array}
 \right)$. That is, for some $k$, $A=\left(
 \begin{array}{cc}
 ka_{12} & a_{12}\\
 ka_{22} & a_{22}
 \end{array}
 \right)$.

Let $P_2=\left(
 \begin{array}{cc}
 1 & 0\\
 -k & 1
 \end{array}
 \right)$. Then $A$ is congruent to $P_2^tAP_2=\left(
 \begin{array}{cc}
 0 & a_{12}-ka_{22}\\
 0 & a_{22}
 \end{array}
 \right)=\left(
 \begin{array}{cc}
 0 & b\\
 0 & a_{22}
 \end{array}
 \right)$.

If $b\neq 0$, then, letting $P_3=\left(
 \begin{array}{cc}
 b^{-1} & -b^{-1}a_{22}\\
 0 & 1
 \end{array}
 \right)$, we have $P_3^tAP_3=\left(
 \begin{array}{cc}
 0 & 1\\
 0 & 0
 \end{array}
 \right)$. Hence we get Type (1).

If $b=0$, then, letting $P_4=\left(
 \begin{array}{cc}
 1 & 0\\
 0 & x
 \end{array}
 \right)$, we have $P_4^tAP_4=\left(
 \begin{array}{cc}
 0 & 0\\
 0 & a_{22}x^2
 \end{array}
 \right)$. Choosing suitable $x$, $a_{22}x^2$ can be $1$ or $\eta$. Hence we get Type (2).

 Let $P=\left(
 \begin{array}{cc}
x_{11} & x_{12}\\
 x_{21} & x_{22}
 \end{array}
 \right)$.
By calculation, the $(2,2)^{th}$ element of $P^t \diag(0,\nu)P$ is
$\nu x_{22}^2$. Hence matrices of Type (2) with different $\nu$ are
not congruent. Since the matrix is symmetric for Type (2), and not
symmetric for Type (1), matrices with different types are not
congruent. \qed

\medskip

\begin{thm}\label{isomorphic-1}
Suppose that $p>2$, $G$ and $\bar{G}$ are finite $p$-groups such that $d(G)=3$, $\Phi(G)\le Z(G)$, $G'\cong
C_p^3$, and the type of $G/G'$ be $(p^{m_1},p^{m_2},p^{m_3})$ where $m_1\ge m_2\ge m_3$. Two characteristic matrices of $G$ and $\bar{G}$ are $w(G)=(w_{ij})$ and $w(\bar{G})=(\bar{w}_{ij})$ respectively.
Then $G\cong \bar{G}$ if and only if there exists $X={\scriptsize\left(
 \begin{array}{ccc}
 x_{11}& x_{12}&x_{13}\\
 x_{21}p^{m_1-m_2}& x_{22}&x_{23}\\
 x_{31}p^{m_1-m_3}& x_{32}p^{m_2-m_3}&x_{33}
 \end{array}
 \right)}$, an invertible matrix over $F_p$, such that $w(\bar{G})=\det(X)^{-1}X_2w(G)X^t$, where
$X_2={\scriptsize\left(
 \begin{array}{ccc}
 x_{11}& x_{12}p^{m_1-m_2}&x_{13}p^{m_1-m_3}\\
 x_{21}& x_{22}&x_{23}p^{m_2-m_3}\\
 x_{31}& x_{32}&x_{33}
 \end{array}
 \right)}$.
 \end{thm}
 \demo
Let $a_1,a_2,a_3\rg$ and $\bar{a}_1,\bar{a}_2,\bar{a}_3$ be two set of characteristic generators about $w(G)$ and $w(\bar{G})$ respectively. Suppose that $\theta$ is an isomorphism from
$\bar{G}$ to $G$. We may assume that
$$\bar{a}_1^\theta\equiv a_1^{x_{11}}a_2^{x_{12}}a_3^{x_{13}}\ (\mod G'),\ \ \bar{a}_2^\theta \equiv
a_1^{x_{21}p^{m_1-m_2}}a_2^{x_{22}}a_3^{x_{23}}\ (\mod G'),$$
$$\bar{a}_3^\theta\equiv
a_1^{x_{31}p^{m_1-m_3}}a_2^{x_{32}p^{m_2-m_3}}a_3^{x_{33}} (\mod
G').$$
By calculation,
\begin{eqnarray*}
  \bar{x}^\theta &=& [\bar{a}_2,\bar{a}_3]^\theta=[\bar{a}_2^\theta,\bar{a}_3^\theta]= [a_1^{x_{21}p^{m_1-m_2}}a_2^{x_{22}}a_3^{x_{23}},a_1^{x_{31}p^{m_1-m_3}}a_2^{x_{32}p^{m_2-m_3}}a_3^{x_{33}}]\\
   &=& x^{x_{22}x_{33}-x_{23}x_{32}p^{m_2-m_3}}y^{-x_{21}x_{33}p^{m_1-m_2}+x_{23}x_{31}p^{m_1-m_3}}z^{x_{21}x_{32}p^{m_1-m_3}-x_{22}x_{31}p^{m_1-m_3}}
\end{eqnarray*}
Similarly we have
\begin{eqnarray*}
\bar{y}^\theta
&=&x^{x_{32}x_{13}p^{m_2-m_3}-x_{33}x_{12}}y^{-x_{31}x_{13}p^{m_1-m_3}+x_{33}x_{11}}z^{x_{31}x_{12}p^{m_1-m_3}-x_{32}x_{11}p^{m_2-m_3}}
\end{eqnarray*}
and
\begin{eqnarray*}\bar{z}^\theta
&=&x^{x_{12}x_{23}-x_{13}x_{22}}y^{-x_{11}x_{23}+x_{13}x_{21}p^{m_1-m_2}}z^{x_{11}x_{22}-x_{12}x_{21}p^{m_1-m_2}}
\end{eqnarray*}
Let 
\begin{eqnarray*}
X_1&=&{\scriptsize\left(
 \begin{array}{ccc}
{x_{22}x_{33}-x_{23}x_{32}p^{m_2-m_3}} &{-x_{21}x_{33}p^{m_1-m_2}+x_{23}x_{31}p^{m_1-m_3}}& {x_{21}x_{32}p^{m_1-m_3}-x_{22}x_{31}p^{m_1-m_3}}\\
{x_{32}x_{13}p^{m_2-m_3}-x_{33}x_{12}} &{-x_{31}x_{13}p^{m_1-m_3}+x_{33}x_{11}}  &{x_{31}x_{12}p^{m_1-m_3}-x_{32}x_{11}p^{m_2-m_3}}\\
 {x_{12}x_{23}-x_{13}x_{22}} &  {-x_{11}x_{23}+x_{13}x_{21}p^{m_1-m_2}}  & {x_{11}x_{22}-x_{12}x_{21}p^{m_1-m_2}}
 \end{array}
 \right)}
 \end{eqnarray*}
Since $(\bar{x}^{\bar{w}_{11}}\bar{y}^{\bar{w}_{12}}\bar{z}^{\bar{w}_{13}})^\theta=(\bar{a}_1^{p^{m_1}})^\theta=
a_1^{x_{11}p^{m_1}}a_2^{x_{12}p^{m_1}}a_3^{x_{13}p^{m_1}}$,
we have
\begin{equation} \label{eq:1}
(\bar{w}_{11},\bar{w}_{12},\bar{w}_{13})X_1=
(x_{11},x_{12}p^{m_1-m_2},x_{13}p^{m_1-m_3})
 {\scriptsize\left(
 \begin{array}{ccc}
 w_{11}& w_{12}&w_{13}\\
 w_{21}& w_{22}&w_{23}\\
 w_{31}& w_{32}&w_{33}
 \end{array}
 \right)}
  \end{equation}
Similarly, by
$(\bar{x}^{\bar{w}_{21}}\bar{y}^{\bar{w}_{22}}\bar{z}^{\bar{w}_{23}})^\theta=(\bar{a}_2^{p^{m_2}})^\theta=
a_1^{x_{21}p^{m_1}}a_2^{x_{22}p^{m_2}}a_3^{x_{23}p^{m_2}}$, we have
\begin{equation} \label{eq:2}
(\bar{w}_{21},\bar{w}_{22},\bar{w}_{23})X_1=
(x_{21},x_{22},x_{23}p^{m_2-m_3})
 {\scriptsize\left(
 \begin{array}{ccc}
 w_{11}& w_{12}&w_{13}\\
 w_{21}& w_{22}&w_{23}\\
 w_{31}& w_{32}&w_{33}
 \end{array}
 \right)}
  \end{equation}
By
$(\bar{x}^{\bar{w}_{31}}\bar{y}^{\bar{w}_{32}}\bar{z}^{\bar{w}_{33}})^\theta=(\bar{a}_3^{p^{m_3}})^\theta=
a_1^{x_{31}p^{m_1}}a_2^{x_{32}p^{m_2}}a_3^{x_{33}p^{m_3}}$, we have
\begin{equation} \label{eq:3}
(\bar{w}_{31},\bar{w}_{32},\bar{w}_{33})X_1=
(x_{31},x_{32},x_{33})
 {\scriptsize\left(
 \begin{array}{ccc}
 w_{11}& w_{12}&w_{13}\\
 w_{21}& w_{22}&w_{23}\\
 w_{31}& w_{32}&w_{33}
 \end{array}
 \right)}
  \end{equation}
Let $X_2={\scriptsize\left(
 \begin{array}{ccc}
 x_{11}& x_{12}p^{m_1-m_2}&x_{13}p^{m_1-m_3}\\
 x_{21}& x_{22}&x_{23}p^{m_2-m_3}\\
 x_{31}& x_{32}&x_{33}
 \end{array}
 \right)}.$ It follow from equations (\ref{eq:1}), (\ref{eq:2}) and (\ref{eq:3}) that
\begin{equation}\label{eq:4}
w(\bar{G})X_1=X_2w(G).
\end{equation}
Let $X={\scriptsize\left(
 \begin{array}{ccc}
 x_{11}& x_{12}&x_{13}\\
 x_{21}p^{m_1-m_2}& x_{22}&x_{23}\\
 x_{31}p^{m_1-m_3}& x_{32}p^{m_2-m_3}&x_{33}
 \end{array}
 \right)}$. Then $X_1=(X^t)^*$, where $X^t$ denote the transpose of $X$ and $X^*$ denote the adjugate of $X$. Right multiplying by $\det(X)^{-1}X^t$, we get
\begin{equation}\label{eq:5}
w(\bar{G})=\det(X)^{-1}X_2w(G)X^t
\end{equation}

On the other hand, if there exists an invertible matrix
$X={\scriptsize\left(
 \begin{array}{ccc}
 x_{11}& x_{12}&x_{13}\\
 x_{21}p^{m_1-m_2}& x_{22}&x_{23}\\
 x_{31}p^{m_1-m_3}& x_{32}p^{m_2-m_3}&x_{33}
 \end{array}
 \right)}$
over $F_p$ satisfying equation (\ref{eq:5}), then
$$\theta:
\bar{a}_1\mapsto a_1^{x_{11}}a_2^{x_{12}}a_3^{x_{13}}, \bar{b}\mapsto
a_1^{x_{21}p^{m_1-m_2}}a_2^{x_{22}}a_3^{x_{23}}, \bar{c}\mapsto
a_1^{x_{31}p^{m_1-m_3}}a_2^{x_{32}p^{m_2-m_3}}a_3^{x_{33}}$$ is an isomorphism from $\bar{G}$ to
$G$.
\qed

\medskip

If $p=2$, then Equation \ref{eq:1}--\ref{eq:3} also hold for $m_2>1$. If $p=2$, $m_1>1$ and $m_2=m_3=1$, then Equation \ref{eq:1} also holds, but Equation \ref{eq:2} and \ref{eq:3} are changed to be
$$
(\bar{w}_{21},\bar{w}_{22},\bar{w}_{23})X_1=
(x_{21},x_{22},x_{23})
 {\scriptsize\left(
 \begin{array}{ccc}
 w_{11}& w_{12}&w_{13}\\
 w_{21}& w_{22}&w_{23}\\
 w_{31}& w_{32}&w_{33}
 \end{array}
 \right)}+(x_{22}x_{23},0,0)
\eqno(2.2')$$ 
and
$$
(\bar{w}_{31},\bar{w}_{32},\bar{w}_{33})X_1=
(x_{31},x_{32},x_{33})
 {\scriptsize\left(
 \begin{array}{ccc}
 w_{11}& w_{12}&w_{13}\\
 w_{21}& w_{22}&w_{23}\\
 w_{31}& w_{32}&w_{33}
 \end{array}
 \right)}+(x_{32}x_{33},0,0)
 \eqno(2.3')$$
 Hence we get the following theorems.

\begin{thm}\label{isomorphic-2}
Suppose that $G$ and $\bar{G}$ are finite $2$-groups such that $d(G)=3$, $\Phi(G)\le Z(G)$, $G'\cong
C_2^3$, and the type of $G/G'$ be $(2^{m_1},2^{m_2},2^{m_3})$ where $m_1\ge m_2\ge m_3$ and $m_2>1$. Two characteristic matrices of $G$ and $\bar{G}$ are $w(G)=(w_{ij})$ and $w(\bar{G})=(\bar{w}_{ij})$ respectively.
Then $G\cong \bar{G}$ if and only if there exists $X={\scriptsize\left(
 \begin{array}{ccc}
 x_{11}& x_{12}&x_{13}\\
 x_{21}2^{m_1-m_2}& x_{22}&x_{23}\\
 x_{31}2^{m_1-m_3}& x_{32}2^{m_2-m_3}&x_{33}
 \end{array}
 \right)}$, an invertible matrix over $F_2$, such that $w(\bar{G})=X_2w(G)X^t$, where
$X_2={\scriptsize\left(
 \begin{array}{ccc}
 x_{11}& x_{12}2^{m_1-m_2}&x_{13}2^{m_1-m_3}\\
 x_{21}& x_{22}&x_{23}2^{m_2-m_3}\\
 x_{31}& x_{32}&x_{33}
 \end{array}
 \right)}$.
 \end{thm}
 
\begin{thm}\label{isomorphic-3}
Suppose that $G$ and $\bar{G}$ are finite $2$-groups such that $d(G)=3$, $\Phi(G)\le Z(G)$, $G'\cong
C_2^3$, and the type of $G/G'$ be $(2^{m_1},2,2)$, where
$m_1>1$. Two characteristic matrix of $G$ and $\bar{G}$ are $w(G)=(w_{ij})$ and $w(\bar{G})=(\bar{w}_{ij})$ respectively.
Then $G\cong \bar{G}$ if and only if there exists $X={\scriptsize\left(
 \begin{array}{ccc}
 1& x_{12}&x_{13}\\
 0 & x_{22}&x_{23}\\
 0 & x_{32}&x_{33}
 \end{array}
 \right)}$, an invertible matrix over $F_2$, such that $w(\bar{G})=X_2w(G)X^t+{\scriptsize\left(
 \begin{array}{ccc}
 0& 0&0\\
 x_{22}x_{23}& 0&0\\
 x_{32}x_{33}&0&0
 \end{array}
 \right)},$ where $X_2={\scriptsize\left(
 \begin{array}{ccc}
 1& 0&0\\
 x_{21}& x_{22}&x_{23}\\
 x_{31}& x_{32}&x_{33}
 \end{array}
 \right)}.$
 \end{thm}

Next we investigate the minimal index and the maximal index of $\mathcal{A}_1$-subgroups of $G$.

\begin{thm}\label{minimal-index}
Suppose that $G$ is a finite $p$-group such that $d(G)=3$, $\Phi(G)\le Z(G)$, $G'\cong
C_p^3$, and the type of $G/G'$ be $(p^{m_1},p^{m_2},p^{m_3})$, where
$m_1\ge m_2\ge m_3$.
Then $m_3\le I_{\min}\le m_3+2$ and

\rr{1} $I_{\min}=m_3$ if and only if there exists a characteristic
matrix $w(G)=(w_{ij})$ of $G$ such that $\rank \left(
 \begin{array}{cc}
 w_{11}& w_{12}\\
 w_{21}& w_{22}
 \end{array}
 \right)=2$;

\rr{2} If $\rank(w(G))\ge 2$, then $m_3\le
I_{\min}\le {m_3+1}$;

\rr{3} If $\rank(w(G))\le 1$, then $m_3+1\le
I_{\min}\le {m_3+2}$, and $I_{\min}=m_3+1$ if and only if there
exists a characteristic matrix $w(G)=(w_{ij})$ of $G$ such that
$\rank \left(
 \begin{array}{cc}
 w_{11}& w_{12}\\
 w_{21}& w_{22}
 \end{array}
 \right)=1$.
\end{thm}
\demo  Let $D$ be an $\mathcal{A}_1$-subgroup of $G$. Since $DG'/G'$
is a subgroup of $G/G'$ and $d(DG'/G')=2$, we have $|G/G':DG'/G'|\ge
p^{m_3}$. Hence
\begin{equation}\label{eq:6}
|G:D|=|G:DG'||DG'/D|=|G/G':DG'/G'||G'/G'\cap D|\ge p^{m_3}.
\end{equation}
On the other hand, $\lg
a_1,a_2\rg$ has index at most $p^{m_3+2}$. Hence $m_3\le I_{\min}\le {m_3+2}$.

(1) If $I_{\min}=m_3$ and $D\in \mathcal{A}_1$ such that
$|G:D|=p^{m_3}$, then, by equation \ref{eq:6}, we have that $G'\le
D$ and $D/G'$ is of type $(p^{m_1},p^{m_2})$. Notice that an element
of maximal order  in  an abelian $p$-group  must be a direct product
factor. We may assume that $D=\lg a_1,a_2\rg$. Since $G'\le D$,
$|\lg a_1^{p^{m_1}},a_2^{p^{m_2}},z\rg|=p^3$. Hence
$\rank{\scriptsize\left(
 \begin{array}{ccc}
 w_{11}& w_{12}&w_{13}\\
 w_{21}& w_{22}&w_{23}\\
 0& 0&1
 \end{array}
 \right)}=3$. Thus $\rank{\left(
 \begin{array}{cc}
 w_{11}& w_{12}\\
 w_{21}& w_{22}
 \end{array}
 \right)}=2$.

 On the other hand, if there exists a $w(G)$ such that $\rank\left(
 \begin{array}{cc}
 w_{11}& w_{12}\\
 w_{21}& w_{22}
 \end{array}
 \right)=2$, then $\lg a_1,a_2\rg$ is of index $p^{m_3}$ and hence $I_{\min}=m_3$.

\medskip

(2) If $I_{\min}=m_3+2$, then $|\lg
a_1^{p^{m_1}},a_2^{p^{m_2}},z\rg|=p$, and hence $\rank{\scriptsize\left(
 \begin{array}{ccc}
 w_{11}& w_{12}&w_{13}\\
 w_{21}& w_{22}&w_{23}\\
 0& 0&1
 \end{array}
 \right)}=1$. It follows that $w(G)={\scriptsize\left(
 \begin{array}{ccc}
 0& 0&w_{13}\\
 0& 0&w_{23}\\
 w_{31}& w_{32}&w_{33}
 \end{array}
\right)}$. Since the rank of $w(G)$ is $\ge 2$, we have
$(w_{13},w_{23})\neq (0,0)$. Without loss of generality, we assume
that $w_{23}\neq 0$. Let $D=\lg a_1a_3, a_2\rg$. It is easy to see that
$D\cap G'\ge \lg x,z\rg$ and hence $D$ is a $\A_1$-subgroup of index
at most $p^{m_3+1}$. This contradicts $I_{\min}=m_3+2$.

\medskip

(3) By (1), $I_{\min}\ge {m_3+1}$. Assume $D$ is a $\mathcal{A}_1$
subgroup of $G$ and $DG'/G'=\lg \bar{g_1}\rg\times\lg \bar{g_2}\rg$,
where $o(\bar{g_1})=p^u$ and $o(\bar{g_2})=p^v$. $\forall g \in D$,
$\exists i,j,k$ such that $g={g_1}^i{g_2}^j[g_1,g_2]^k$. It is easy
to see that $g\in G'$ if and only if $p^u|i$ and $p^v|j$. Hence
$D\cap G'=\lg {g_1}^{p^u},{g_2}^{p^v},[g_1,g_2]\rg$. Notice that
$\rank(w(G))\le 1$, we have $|D\cap G'|\le p^2$.

If $I_{\min}=m_3+1$ and $D\in \mathcal{A}_1$ such that
$|G:D|=p^{m_3+1}$, then, by equation \ref{eq:6}, we have $DG'/G'$ is
of type $(p^{m_1},p^{m_2})$. Hence we may assume that $D=\lg
a_1,a_2\rg$. Since $|\lg a_1^{p^{m_1}},a_2^{p^{m_2}},z\rg|=p^2$,
$\rank{\scriptsize\left(
 \begin{array}{ccc}
 w_{11}& w_{12}&w_{13}\\
 w_{21}& w_{22}&w_{23}\\
 0& 0&1
 \end{array}
 \right)}=2$. Thus $\rank{\left(
 \begin{array}{cc}
 w_{11}& w_{12}\\
 w_{21}& w_{22}
 \end{array}
 \right)}=1$.

 On the other hand, if there is a $w(G)$ such that $\rank\left(
 \begin{array}{cc}
 w_{11}& w_{12}\\
 w_{21}& w_{22}
 \end{array}
 \right)=1$, then $\lg a_1,a_2\rg$ is of index $p^{m_3+1}$ and hence $I_{\min}=m_3+1$.
\qed

\begin{thm}\label{maximal-index}
Suppose that $G$ is a finite $p$-group such that $d(G)=3$, $\Phi(G)\le Z(G)$, $G'\cong
C_p^3$, and the type of $G/G'$ be $(p^{m_1},p^{m_2},p^{m_3})$, where
$m_1\ge m_2\ge m_3$.
Then $m_1\le I_{\max}\le m_1+2$ and

\rr{1} $I_{\max}=m_1+2$ if and only if there exists a characteristic
matrix $w(G)=(w_{ij})$ of $G$ such that $\left(
 \begin{array}{cc}
 w_{22}& w_{23}\\
 w_{32}& w_{33}
 \end{array}
 \right)=0$;

\rr{2} If $\rank(w(G))\le 1$, then $m_1+1\le
I_{\max}\le {m_1+2}$;

\rr{3} If $\rank(w(G))=3$, then $m_1\le I_{\max}\le {m_1+1}$;

\rr{4} If $I_{\max}\neq m_1+2$, then $I_{\max}= m_1+1$ if and only if one of the following conditions holds:

\ \ \ \ \ \rr{a} There is a characteristic matrix $w(G)$ such that $\rank\left(
 \begin{array}{cc}
 w_{22}& w_{23}\\
 w_{32}& w_{33}
 \end{array}
 \right)=1$;

\ \ \ \ \ \  \rr{b} $m_1=m_2+1$ and there is a characteristic matrix $w(G)$ such that $\left(
 \begin{array}{cc}
 w_{11}& w_{13}\\
 w_{31}& w_{33}
 \end{array}
 \right)=0$.
\end{thm}
\demo Assume  $D$ is an $\mathcal{A}_1$-subgroups of $G$. Then
$D\cap \Phi(G)=D\cap Z(G)\le Z(D)=\Phi(D)$. Hence
$\Phi(D)=D\cap\Phi(G)$. Thus $D\Phi(G)/\Phi(G)\cong D/D\cap
\Phi(G)=D/\Phi(D)=p^2$. It follows that $|G/G':DG'/G'|\le p^{m_1}$.
Since $|G'\cap D|\ge p$, we have
\begin{equation}\label{eq:7}
|G:D|=|G:DG'||DG'/D|=|G/G':DG'/G'||G'/G'\cap D|\le p^{m_1+2}.
\end{equation}
On the other hand, $\lg
a_2,a_3\rg$ has index at least $p^{m_1}$. Hence $m_1\le I_{\max}\le {m_1+2}$.

\smallskip
(1) If $I_{\max}=m_1+2$ and $D\in \mathcal{A}_1$ such that
$|G:D|=p^{m_1+2}$, then, by equation \ref{eq:7}, we have that
$|G'\cap D|=p$ and $DG'/G'$ is of type $(p^{m_2},p^{m_3})$. Notice
that an element of minimal order  in  an abelian $p$-group  must be
a direct product factor if it is a generator. We may assume that
$D=\lg a_2,a_3\rg$. Since $|\lg
x,a_2^{p^{m_2}},a_3^{p^{m_3}}\rg|=p$, $\rank{\scriptsize\left(
 \begin{array}{ccc}
 1&0&0\\
 w_{21}& w_{22}&w_{23}\\
 w_{31}& w_{32}&w_{33}
 \end{array}
 \right)}=1$. Hence ${\left(
 \begin{array}{cc}
 w_{22}& w_{23}\\
 w_{32}& w_{33}
 \end{array}
 \right)}=0$.

 On the other hand, if there is a $w(G)$ such that $\left(
 \begin{array}{cc}
 w_{22}& w_{23}\\
 w_{32}& w_{33}
 \end{array}
 \right)=0$, then $\lg a_2,a_3\rg$ is of index $p^{m_1+2}$ and hence $I_{\max}=m_1+2$.

\smallskip

(2) Assume the contrary, then $I_{\max}=m_1$. In this case, $|G:\lg a_2,a_3\rg|=p^{m_1}$. It follows that $|\lg x,a_2^{p^{m_2}},a_3^{p^{m_3}}\rg|=p^3$. Then the matrix ${\scriptsize\left(
 \begin{array}{ccc}
 1&0&0\\
 w_{21}& w_{22}&w_{23}\\
 w_{31}& w_{32}&w_{33}
 \end{array}
 \right)}$ is invertible. Hence $\rank(w(G))\ge 2$, a contradiction. So $I_{\max}\ge m_1+1$.

\medskip

(3) By (1), if $I_{\max}=m_1+2$, then $\rank(w(G))\le 2$.

\medskip

(4) If (a) holds, then $|\lg x,a_2^{p^{m_2}},a_3^{p^{m_3}}\rg|=p^2$.
Hence $\lg a_2,a_3\rg$ is of index $p^{m_1+1}$. If (b) holds, then
$|\lg y,a_1^{p^{m_1}},a_2^{p^{m_2}}\rg|=p$. Hence $\lg a_1,a_3\rg$
is of index $p^{m_2+2}=p^{m_1+1}$. So $I_{\max}=m_1+1$.

\medskip

 On the other hand, if $I_{\max}=m_1+1$ and
 $D\in \mathcal{A}_1$ such that $|G:D|=p^{m_1+1}$, then, $|D|=p^{m_2+m_3+2}$. Since $|D/G'\cap D|=|DG'/G'|\ge p^{m_2+m_3}$, $|G'\cap D|\le p^2$.

\medskip

 Case (a): There exists a $D\in\mathcal{A}_1$ such that $|G:D|=p^{m_1+1}$ and $|G'\cap D|=p^2$.

  Since $|DG'/G'|=|D/G'\cap D|=p^{m_2+m_3}$, $DG'/G'$ is of type $(p^{m_2},p^{m_3})$. Hence we may assume that $D=\lg a_2,a_3\rg$. Since $|G'\cap D|=p^2$, $\rank{\scriptsize\left(
 \begin{array}{ccc}
 1&0&0\\
 w_{21}& w_{22}&w_{23}\\
 w_{31}& w_{32}&w_{33}
 \end{array}
 \right)}=2$. Hence $\rank\left(
 \begin{array}{cc}
 w_{22}& w_{23}\\
 w_{32}& w_{33}
 \end{array}
 \right)=1$. Condition (a) holds.

  \medskip

Case (b): For all $D\in\mathcal{A}_1$ such that $|G:D|=p^{m_1+1}$, we have $|G'\cap D|=p$.

Then $|DG'/G'|=|D/G'\cap D|=p^{m_2+m_3+1}$. It follows that $DG'/G'$ is of type $(p^{m_2+1},p^{m_3})$ or $(p^{m_2},p^{m_3+1})$.

Subcase (b1): $DG'/G'$ is of type $(p^{m_2+1},p^{m_3})$.

If $G=\lg D,a_1\rg$, then we may assume that $D=\lg a_2b,a_3\rg$. Since $|\lg x,a_3^{p^{m_3}}\rg|=p$, we have $|\lg x,a_2^{p^{m_2}},a_3^{p^{m_3}}\rg|\le p^2$. That is, $|G'\cap \lg a_2,a_3\rg|\le p^2$. Hence $|G:\lg a_2,a_3\rg|\ge p^{m_1+1}$. Since $I_{\max}=m_1+1$, we have $|G:\lg a_2,a_3\rg|= p^{m_1+1}$ and $|G'\cap \lg a_2,a_3\rg|= p^2$, which contradict the hypothesis of Case (b).

If $G=\lg D,a_2\rg$, then we may assume that $m_1>m_2$. (If $m_1=m_2$, it is reduced to the case $G=\lg D,a_1\rg$.) Since $|G/G':DG'/G'|\le p^{m_2}$, we get $|G:D|=|G:DG'||DG':D|\le p^{m_2}|G'/G'\cap D|=p^{m_2+2}$. Since $I_{\max}=m_1+1$, we have $m_1=m_2+1$ and $|G:DG'|=p^{m_2}$.
Since $DG'/G'$ is of type $(p^{m_1},p^{m_3})$, we may assume that $D=\lg a_1, a_3\rg$. Since $|G'\cap D|=p$, we have $\rank{\scriptsize\left(
 \begin{array}{ccc}
 w_{11}& w_{12}&w_{13}\\
 0& 1 & 0\\
 w_{31}& w_{32}&w_{33}
 \end{array}
 \right)}=1$. Hence $\left(
 \begin{array}{cc}
 w_{11}& w_{13}\\
 w_{31}& w_{33}
 \end{array}
 \right)=0$. Condition (b) holds.

If $G=\lg D,a_3\rg$, then we may assume that $D=\lg a_1,a_2\rg$. Hence $m_1=m_2+1$ and $m_2=m_3$. It is reduced to the case $G=\lg D,a_2\rg$.

\medskip

Subcase (b2): $DG'/G'$ is of type $(p^{m_2},p^{m_3+1})$.

If $m_2=m_3$, then it can be reduced to Subcase (b1). Hence we may
assume that $m_2>m_3$.

If $G=\lg D,a_1\rg$, then we may assume that
$D=\lg a_2,a_3c\rg$. Since $|\lg x,a_2^{p^{m_2}}\rg|=p$, we have
that $|\lg x,a_2^{p^{m_2}},a_3^{p^{m_3}}\rg|\le p^2$. That is,
$|G'\cap \lg a_2,a_3\rg|\le p^2$. Hence $|G:\lg a_2,a_3\rg|\ge
p^{m_1+1}$. Since $I_{\max}=m_1+1$, we have $|G:\lg a_2,a_3\rg|=
p^{m_1+1}$ and $|G'\cap \lg a_2,a_3\rg|= p^2$, which contradict the
hypothesis of Case (b).

If $G=\lg D,a_2\rg$, then we may assume $D=\lg a_1, a_3c\rg$. It
follows that $m_1=m_2$, and it can be reduced to the case $G=\lg
D,a_1\rg$.

If $G=\lg D,a_3\rg$, then we may assume $D=\lg a_1, a_2\rg$. It
follows that $m_1=m_2=m_3+1$. Hence $I_{\max}={m_3+2}$. Notice that
$|G/G':DG'/G'|=p^{m_3}$, we have $|D\cap G'|=p^2$, which also contradict the
hypothesis of Case (b). \qed
 \section{The case $m_1>m_2>m_3$}

In this section, we deal with the case $m_1>m_2>m_3$. In this case, $p$ is either odd or $2$.

\begin{thm}\label{1>2>3} Let $G$ be a finite $p$-group with $d(G)=3$, $\Phi(G)\le Z(G)$ and $G'\cong C_p^3$. If the type of $G/G'$ is $(p^{m_1},p^{m_2},p^{m_3})$ where $m_1>m_2>m_3$, then we can choose suitable generators of $G$, such that the characteristic matrix of $G$ is one of the
following matrices, and different matrices give non-isomorphic groups. $($Where $\eta$ is a fixed square
non-residue modulo odd $p$, $\nu=1$ or $\eta$, and $t\neq 0. )$

\medskip

{\rm (A1)} ${\scriptsize\left(
 \begin{array}{ccc}
 1& 0&0\\
 0& \nu_1& 0\\
 0& 0& \nu_2
 \end{array}
 \right)}$ where $\nu_1,\nu_2=1$ or $\eta$\ \
 {\rm (A2)} ${\scriptsize\left(
 \begin{array}{ccc}
 1& 0&0\\
 0& 0& 1\\
 0& t& 0
 \end{array}
 \right)}$\ \
{\rm (A3)} ${\scriptsize\left(
 \begin{array}{ccc}
 0& 0&  t\\
 0& 1&0\\
 1 & 0&0
 \end{array}
 \right)}$\ \

 \smallskip

{\rm (A4)} ${\scriptsize\left(
 \begin{array}{ccc}
 0& 0&1\\
 1& 0&0\\
 0& 1&0
 \end{array}
 \right)}$\
  {\rm (A5)} ${\scriptsize\left(
 \begin{array}{ccc}
 0& 1&0\\
 0& 0 &1\\
 1& 0&0
 \end{array}
 \right)}$\
 {\rm (A6)}  ${\scriptsize\left(
 \begin{array}{ccc}
 0& t&0\\
 1& 0&0\\
 0& 0&1
 \end{array}
 \right)}$\ \

\smallskip

{\rm (B1)} ${\scriptsize\left(
 \begin{array}{ccc}
 1& 0&0\\
 0& \nu& 0\\
 0& 0& 0
 \end{array}
 \right)}$\ \
 {\rm (B2)}\ ${\scriptsize\left(
 \begin{array}{ccc}
 1& 0&  0\\
 0& 0&1\\
 0 & 0&0
 \end{array}
 \right)}$\ \
 {\rm (B3)} ${\scriptsize\left(
 \begin{array}{ccc}
 0& 0&1\\
 0& 1& 0\\
 0& 0& 0
 \end{array}
 \right)}$\ \
  {\rm (B4)} ${\scriptsize\left(
 \begin{array}{ccc}
 0& 1&0\\
 0& 0& 1\\
 0& 0& 0
 \end{array}
 \right)}$\ \
 {\rm (B5)} ${\scriptsize\left(
 \begin{array}{ccc}
 0& 0& 1\\
 0& 0& 0\\
 0& 1& 0
 \end{array}
 \right)}$
 \smallskip

{\rm (B6)} ${\scriptsize\left(
 \begin{array}{ccc}
 1& 0&0\\
 0& 0& 0\\
 0& 0& \nu
 \end{array}
 \right)}$\ \
{\rm (B7)} ${\scriptsize\left(
 \begin{array}{ccc}
 1& 0&0\\
 0& 0&0\\
 0& 1&0
 \end{array}
 \right)}$\ \
  {\rm (B8)}  ${\scriptsize\left(
 \begin{array}{ccc}
 0& 1&0\\
 0& 0&0\\
 0& 0&1
 \end{array}
 \right)}$\ \
{\rm (B9)} ${\scriptsize\left(
 \begin{array}{ccc}
 0& 0&  0\\
 0& 1&0\\
 1 & 0&0
 \end{array}
 \right)}$
{\rm (B10)} ${\scriptsize\left(
 \begin{array}{ccc}
 0& 0&0\\
 0& 0&1\\
1& 0&0
 \end{array}
 \right)}$

 \smallskip

   {\rm (B11)} ${\scriptsize\left(
 \begin{array}{ccc}
 0& 0& 0\\
 1& 0& 0\\
 0& 0& 1
 \end{array}
 \right)}$\ \
 {\rm (B12)} ${\scriptsize\left(
 \begin{array}{ccc}
 0& 0& 0\\
 1& 0& 0\\
 0& 1& 0
 \end{array}
 \right)}$\ \
{\rm (B13)} ${\scriptsize\left(
 \begin{array}{ccc}
 0& 1& 0\\
 0& 0& 0\\
 1& 0& 0
 \end{array}
 \right)}$\ \
{\rm (B14)} ${\scriptsize\left(
 \begin{array}{ccc}
 0& 0&t\\
 0& 0 &0\\
 1& 0&0
 \end{array}
 \right)}$\ \

\smallskip

{\rm (B15)} ${\scriptsize\left(
 \begin{array}{ccc}
 0& 0&1\\
 1& 0& 0\\
 0& 0& 0
 \end{array}
 \right)}$\ \
{\rm (B16)} ${\scriptsize\left(
 \begin{array}{ccc}
 0& t& 0\\
 1& 0& 0\\
 0& 0& 0
 \end{array}
 \right)}$\ \
 {\rm (B17)} ${\scriptsize\left(
 \begin{array}{ccc}
 0& 0&0\\
 0& 1 &0\\
 0& 0&\nu
 \end{array}
 \right)}$\ \
 {\rm (B18)}  ${\scriptsize\left(
 \begin{array}{ccc}
 0& 0&0\\
 0& 0&1\\
 0& t&0
 \end{array}
 \right)}$\ \

 \smallskip

 {\rm (C1)}  ${\scriptsize\left(
 \begin{array}{ccc}
 1& 0&0\\
 0& 0&0\\
 0& 0&0
 \end{array}
 \right)}$\ \
 {\rm (C2)}  ${\scriptsize\left(
 \begin{array}{ccc}
 0& 1&0\\
 0& 0&0\\
 0& 0&0
 \end{array}
 \right)}$\ \
 {\rm (C3)} ${\scriptsize\left(
 \begin{array}{ccc}
 0& 0& 1\\
 0& 0& 0\\
 0& 0& 0
 \end{array}
 \right)}$\ \
 {\rm (C4)} ${\scriptsize\left(
 \begin{array}{ccc}
 0& 0& 0\\
 1& 0& 0\\
 0& 0& 0
 \end{array}
 \right)}$\ \
 {\rm (C5)} ${\scriptsize\left(
 \begin{array}{ccc}
 0& 0&0\\
 0& 1& 0\\
 0& 0& 0
 \end{array}
 \right)}$\ \

 \smallskip

 {\rm (C6)} ${\scriptsize\left(
 \begin{array}{ccc}
 0& 0& 0\\
 0& 0& 1\\
 0& 0& 0
 \end{array}
 \right)}$\ \
 {\rm (C7)} ${\scriptsize\left(
 \begin{array}{ccc}
 0& 0& 0\\
 0& 0& 0\\
 0& 0& 1
 \end{array}
 \right)}$\ \
 {\rm (C8)} ${\scriptsize\left(
 \begin{array}{ccc}
 0& 0& 0\\
 0& 0& 0\\
 0& 1& 0
 \end{array}
 \right)}$\ \
{\rm (C9)} ${\scriptsize\left(
 \begin{array}{ccc}
 0& 0&0\\
 0& 0& 0\\
 1& 0& 0
 \end{array}
 \right)}$
 {\rm (C10)} ${\scriptsize\left(
 \begin{array}{ccc}
 0& 0& 0\\
 0& 0& 0\\
 0& 0& 0
 \end{array}
 \right)}$\ \
\end{thm}
\demo
Suppose that $G$ and $\bar{G}$ are two groups described in the theorem. By Theorem \ref{isomorphic-1} and \ref{isomorphic-2}, $\bar{G}\cong G$ if and only if there exist invertible matrices $X_2={\scriptsize\left(
 \begin{array}{ccc}
 x_{11}&      0&      0\\
 x_{21}& x_{22}&      0\\
 x_{31}& x_{32}&   x_{33}
 \end{array}
 \right)}$ and $X={\scriptsize\left(
 \begin{array}{ccc}
x_{11} & x_{12}  & x_{13}\\
0 & x_{22}  &x_{23}\\
0 &  0  & x_{33}
 \end{array}
 \right)}$ such that
$w(\bar{G})=\det(X)^{-1}X_2w(G)X^t.$

Notice that $X_2={\scriptsize\left(
 \begin{array}{ccc}
 x_{11}&      0&      0\\
 0    & x_{22}&      0\\
 0    & 0     &   x_{33}
 \end{array}
 \right)}{\scriptsize\left(
 \begin{array}{ccc}
 1&      0&      0\\
 x_{21}x_{22}^{-1}& 1&      0\\
 x_{31}x_{33}^{-1} &x_{32}x_{33}^{-1}&1
 \end{array}
 \right)}$, and $X^t$ has same decomposition. We can use following three kinds of operation to simplify $w(G)$.

Operation I. Take  $X_2=X=\diag(x_{11},x_{22},x_{33})$.  Then $w(\bar{G})=\det(X)^{-1}Xw(G)X^t$.

Operation II. Take $X_2={\scriptsize\left(
 \begin{array}{ccc}
 1&      0&      0\\
 x_{21}& 1&      0\\
 x_{31} &x_{32}&1
 \end{array}
 \right)}$ and $X=E_3$. Then $w(\bar{G})=X_2w(G)$. This operation means adding the 1st row multiplied by a scalar $x_{21}$ to the 2nd row, adding the 1st row multiplied by a scalar $x_{31}$ to the 3rd row, and adding the 2nd row multiplied by a scalar $x_{32}$ to the 3rd row.

Operation III. Take $X_2=E_3$ and $X^t={\scriptsize\left(
 \begin{array}{ccc}
 1&      0&      0\\
 x_{12}& 1&      0\\
 x_{13} &x_{23}& 1
 \end{array}
 \right)}$. Then $w(\bar{G})=w(G)X^t$. This operation means adding the 3rd column multiplied by a scalar $x_{13}$ to the 1st column, adding the 3rd column multiplied by a scalar $x_{23}$ to the 2nd column, and adding the 2nd column multiplied by a scalar $x_{12}$ to the 1st row.

Case 1. $w(G)$ is invertible.

In this case, the 3rd column of $w(G)$ is not $(0,0,0)^t$. Hence we have three types of $w(G)$:

{\rm (a)} ${\scriptsize\left(
 \begin{array}{ccc}
 *& *&w_{13}\\
 *& *&*\\
 *& *&*
 \end{array}
 \right)}$ where $w_{13}\neq 0$;\ \
  {\rm (b)} ${\scriptsize\left(
 \begin{array}{ccc}
 *& *&0\\
 *& * &w_{23}\\
 *& *&*
 \end{array}
 \right)}$ where $w_{23}\neq 0$;

 {\rm (c)} ${\scriptsize\left(
 \begin{array}{ccc}
 *& *&0\\
 *& * &0\\
 *& *&w_{33}
 \end{array}
 \right)}$ where $w_{33}\neq 0$.

Using operation I, II and III, we can not change the type of $w(G)$. Hence matrices of different types determine non-isomorphic groups.

Now we assume that $w(G)$ is of Type (b). (All types are almost the same.)
Using operation II and III, $w(G)$ can be simplified to be ${\scriptsize\left(
 \begin{array}{ccc}
 *& *&0\\
 0& 0 &w_{23}\\
 *& *&0
 \end{array}
 \right)}$. Since the 2nd column is not $(0,0,0)^t$, there are two types of $w(G)$:

  {\rm (b1)} ${\scriptsize\left(
 \begin{array}{ccc}
 *& w_{12}&0\\
 0& 0 &w_{23}\\
 *& *&0
 \end{array}
 \right)}$ where $w_{12}\neq 0$;\ \ {\rm (b2)} ${\scriptsize\left(
 \begin{array}{ccc}
 *& 0&0\\
 0& 0 &w_{23}\\
 *& w_{32}&0
 \end{array}
 \right)}$ where $w_{32}\neq 0$.

Using operation I, II and III, we can not change the type of $w(G)$ from (b1) to (b2). Hence matrices of different types determine non-isomorphic groups.

Using operation II and III, $w(G)$ can be simplified to be \rr{b$1'$} ${\scriptsize\left(
 \begin{array}{ccc}
 0& w_{12}&0\\
 0& 0 &w_{23}\\
 w_{31}& 0&0
 \end{array}
 \right)}$ or \rr{b$2'$} ${\scriptsize\left(
 \begin{array}{ccc}
 w_{11}& 0&0\\
 0& 0& w_{23}\\
 0& w_{32}& 0
 \end{array}
 \right)}$.

If $w(G)$ is the matrix of Type (b$1'$), then, using operation I, where $X=\diag(w_{23},w_{31},w_{12})$, we have $w(\bar{G})=\det(X)^{-1}Xw(G)X^t={\scriptsize\left(
 \begin{array}{ccc}
 0& 1&0\\
 0& 0 &1\\
 1& 0&0
 \end{array}
 \right)}$. Hence we get the  matrix (A5).

 If $w(G)$ is the matrix of Type (b$2'$), then, using operation I, where $X=\diag(w_{23},w_{23},w_{11})$, we have $w(\bar{G})=\det(X)^{-1}Xw(G)X^t={\scriptsize\left(
 \begin{array}{ccc}
 1& 0& 0\\
 0& 0 &1\\
 0&  w_{32}w_{23}^{-1} &0
 \end{array}
 \right)}$. Hence we get the matrix (A2).

Obviously, two groups with characteristic matrices $w(G)$ and $w(\bar{G})$ of Type (A2) are mutually isomorphic if and only if there exsits an invertible matrix $X=\diag(x_{11},x_{22},x_{33})$ such that $w(\bar{G})=\det(X)^{-1}Xw(G)X^t$.
If ${\scriptsize\left(
 \begin{array}{ccc}
 1& 0& 0\\
 0& 0 &1\\
 0& \bar{t} &0
 \end{array}
 \right)}=\det(X)^{-1}X{\scriptsize\left(
 \begin{array}{ccc}
 1& 0& 0\\
 0& 0 &1\\
 0& t &0
 \end{array}
 \right)}X^t$, then, by calculation, we have $x_{11}=1$ and $\bar{t}=t$. Hence matrices with different $t$ give non-isomorphic groups.

Similar arguments as above give characteristic matrices (A1)--(A6).

\medskip

Case 2. $w(G)$ is not invertible.

Similar to Case 1, using Operation II and III, $w(G)$ can be simplified to be such a matrix, in which every column and every row have at most one non-zero element.
If $\rank(w(G))=2$, then we will get characteristic matrices (B1)--(B18). If $\rank(w(G))=1$, then we will get characteristic matrices (C1)--(C9). $w(G)=0$ is the matrix (C10). Also similar to Case 1, we can prove that different matrices give non-isomorphic groups.
\qed

\medskip

Next we give an example to explain how to use Theorem \ref{minimal-index} and \ref{maximal-index} to get $I_{\min}$ and $I_{\max}$.

\begin{thm}\label{B18-property} Let $G$ be a finite $p$-group determined by the characteristic matrix {\rm (B18)} in Theorem $\ref{1>2>3}.$
Then $I_{\min}=m_3+1$, $I_{\max}=m_1$ for $m_1>m_2+1$, and $I_{\max}=m_1+1$ for $m_1=m_2+1$.
\end{thm}
\demo
Since $\rank(w(G))=2$, by Theorem \ref{minimal-index} (2), $m_3\le I_{\min} \le m_3+1$. By Theorem \ref{minimal-index} (1), $I_{\min} =m_3$ if and only if there exist invertible matrices $X_2={\scriptsize\left(
 \begin{array}{ccc}
 x_{11}&      0&      0\\
 x_{21}& x_{22}&      0\\
 x_{31}& x_{32}&   x_{33}
 \end{array}
 \right)}$ and $X={\scriptsize\left(
 \begin{array}{ccc}
x_{11} & x_{12}  & x_{13}\\
0 & x_{22}  &x_{23}\\
0 &  0  & x_{33}
 \end{array}
 \right)}$ and
$w(\bar{G})=\det(X)^{-1}X_2w(G)X^t=(\bar{w}_{ij})$ such that $\left(
 \begin{array}{cc}
 \bar{w}_{11}& \bar{w}_{12}\\
 \bar{w}_{21} &\bar{w}_{22}
 \end{array}
 \right)$ is invertible. By calculation, we have
 \begin{equation*}\label{B18}
 \det(X)w(\bar{G})={\scriptsize\left(
 \begin{array}{ccc}
 0&      0&      0\\
 x_{22}x_{13}& x_{22}x_{23}&  x_{22}x_{33}\\
 tx_{33}x_{12}+x_{32}x_{13}&tx_{33}x_{22}+x_{32}x_{23}&   x_{32}x_{33}
 \end{array}
 \right)}
 \end{equation*}
 Since $\left(
 \begin{array}{cc}
 \bar{w}_{11}& \bar{w}_{12}\\
 \bar{w}_{21} &\bar{w}_{22}
 \end{array}
 \right)$ is not invertible, we have $I_{\min}=m_3+1$.

Since $\left(
 \begin{array}{cc}
 \bar{w}_{22}& \bar{w}_{23}\\
 \bar{w}_{32} &\bar{w}_{33}
 \end{array}
 \right)$ is invertible, by Theorem \ref{maximal-index}, $I_{\max}\neq m_1+2$ and Condition (a) in Theorem \ref{maximal-index} (4) does not hold. If we take $x_{12}=x_{32}=0$, then $\left(
 \begin{array}{cc}
 \bar{w}_{11}& \bar{w}_{13}\\
 \bar{w}_{31} &\bar{w}_{33}
 \end{array}
 \right)=0$. Hence Condition (b) in Theorem \ref{maximal-index} (4) holds if and only if $m_1=m_2+1$. Hence $I_{\max}=m_1$ for $m_1>m_2+1$, and $I_{\max}=m_1+1$ for $m_1=m_2+1$.\qed

\medskip

Similar method as Theorem \ref{B18-property} gives the following theorem. The details are omitted.

\begin{thm}\label{1>2>3-property} Let $G$ be a finite $p$-group determined by a characteristic matrix in Theorem $\ref{1>2>3}.$

\rr{1} If $I_{\min}={m_3}$, then the characteristic matrix is one of the following: {\rm (A1)--(A6), (B1)--(B4), (B15)--(B16)};

\rr{2} If $I_{\min}={m_3+1}$, then the characteristic matrix is one of the following: {\rm (B5)--(B14), (B17)--(B18), (C1)--(C6)};

\rr{3} If $I_{\min}={m_3+2}$, then the characteristic matrix is one of the following: {\rm (C7)--(C10)};

\rr{4} If $I_{\max}={m_1+2}$, then the characteristic matrix is one of the following: {\rm (B13)--(B16), (C1)--(C4), (C9)--(C10)};

\rr{5} If $I_{\max}={m_1+1}$, then the characteristic matrix is one of the following: {\rm (A3)--(A6), (B1)--(B12), (B18) for $m_1=m_2+1$, (C5)--(C8)};

\rr{6} If $I_{\max}={m_1}$, then the characteristic matrix is one of the following: {\rm (A1)--(A2), (B17), (B18) for $m_1>m_2+1$}.
\end{thm}

 \section{The case $m_1>m_2=m_3$}

In this case, we assume that $m_2>1$ for $p=2$.

\begin{thm}\label{1>2=3} Let $G$ be a finite $p$-group with $d(G)=3$, $\Phi(G)\le Z(G)$ and $G'\cong C_p^3$. If the type of $G/G'$ is $(p^{m_1},p^{m_2},p^{m_3})$ where $m_2>1$ for $p=2$ and $m_1>m_2=m_3$, then we can choose suitable generators of $G$, such that the characteristic matrix of $G$ is one of the
following matrices, and different matrices give non-isomorphic groups. Where $\eta$ is a fixed square
non-residue modulo odd $p$, $\nu=1$ or $\eta$, $t\neq 0$, and  $r=1,2,\dots,p-2$.

\smallskip

\rr{1} For odd $p$

\medskip

\rr{D1} ${\scriptsize\left(
 \begin{array}{ccc}
 1 & 0 &0\\
 0&  0& 1\\
 0& -1& 0
 \end{array}
 \right)}$\ \
{\rm (D2)} ${\scriptsize\left(
 \begin{array}{ccc}
1&0&0\\
0& \nu& 1\\
 0& -1& 0
 \end{array}
 \right)}$\ \
{\rm (D3)} ${\scriptsize\left(
 \begin{array}{ccc}
 1  & 0  & 0\\
 0 &1& 0\\
 0 &0 & \nu
 \end{array}
 \right)}$\ \
{\rm (D4)} ${\scriptsize\left(
 \begin{array}{ccc}
 1  & 0  & 0\\
 0 &1& 1\\
 0 &-1& r
 \end{array}
 \right)}$\ \

\smallskip

\rr{E1} ${\scriptsize\left(
 \begin{array}{ccc}
 0 & 0 &0\\
 0&  0& 1\\
 0& -1& 0
 \end{array}
 \right)}$\ \
{\rm (E2)} ${\scriptsize\left(
 \begin{array}{ccc}
0&0&0\\
0& 1& 1\\
 0& -1& 0
 \end{array}
 \right)}$\ \
{\rm (E3)} ${\scriptsize\left(
 \begin{array}{ccc}
 0  & 0  & 0\\
 0 &1& 0\\
 0 &0 & \nu
 \end{array}
 \right)}$\ \
{\rm (E4)} ${\scriptsize\left(
 \begin{array}{ccc}
 0  & 0  & 0\\
 0 &1& 1\\
 0 &-1& r
 \end{array}
 \right)}$

\medskip

\rr{2} For $p=2$

\medskip

{\rm (D5)} ${\scriptsize\left(
 \begin{array}{ccc}
1&0&0\\
0& 0& 1\\
 0& 1& 0
 \end{array}
 \right)}$\ \
\rr{D6} ${\scriptsize\left(
 \begin{array}{ccc}
 1 & 0 &0\\
 0&  1& 0\\
 0& 0& 1
 \end{array}
 \right)}$\ \
{\rm (D7)} ${\scriptsize\left(
 \begin{array}{ccc}
 1  & 0  & 0\\
 0 &1& 0\\
 0 &1 & 1
 \end{array}
 \right)}$

\smallskip

{\rm (E5)} ${\scriptsize\left(
 \begin{array}{ccc}
0&0&0\\
0& 0& 1\\
 0& 1& 0
 \end{array}
 \right)}$\ \
 \rr{E6} ${\scriptsize\left(
 \begin{array}{ccc}
 0 & 0 &0\\
 0&  1& 0\\
 0& 0& 1
 \end{array}
 \right)}$\ \
{\rm (E7)} ${\scriptsize\left(
 \begin{array}{ccc}
 0  & 0  & 0\\
 0 &1& 0\\
 0 &1 & 1
 \end{array}
 \right)}$

 \medskip

 \rr{3} For both odd $p$ and $p=2$

 \medskip

\rr{D8} ${\scriptsize\left(
 \begin{array}{ccc}
 0&   0&   1\\
 0&  1&   0\\
 t & 0 & 0
 \end{array}
 \right)}$\ \
 \rr{D9} ${\scriptsize\left(
 \begin{array}{ccc}
 0&   0&   1\\
 1&  0&   0\\
 0 & 1& 0
 \end{array}
 \right)}$\ \

\smallskip

{\rm (E8)} ${\scriptsize\left(
 \begin{array}{ccc}
 1 &0 &0\\
 0&0& 1\\
 0&0& 0
 \end{array}
 \right)}$\ \
{\rm (E9)} ${\scriptsize\left(
 \begin{array}{ccc}
 1&0&0\\
 0& 0& 0\\
 0&0& \nu
 \end{array}
 \right)}$\ \
 \rr{E10} ${\scriptsize\left(
 \begin{array}{ccc}
 0&   0&    0\\
 0&  0&   1\\
 1 & 0 & 0
 \end{array}
 \right)}$\ \
\rr{E11} ${\scriptsize\left(
 \begin{array}{ccc}
 0&   0&    0\\
 1&  0&   0\\
 0 & 0 & 1
 \end{array}
 \right)}$\ \

 \smallskip

\rr{E12} ${\scriptsize\left(
 \begin{array}{ccc}
 0&   0&   1\\
 0&  1&   0\\
 0 & 0 & 0
 \end{array}
 \right)}$\ \
 \rr{E13} ${\scriptsize\left(
 \begin{array}{ccc}
 0&   0&   1\\
 0&  0&   0\\
 0 & 1& 0
 \end{array}
 \right)}$\ \
 \rr{E14} ${\scriptsize\left(
 \begin{array}{ccc}
 0&   0&   1\\
 1&  0&   0\\
 0 & 0& 0
 \end{array}
 \right)}$\ \
\rr{E15} ${\scriptsize\left(
 \begin{array}{ccc}
 0&   0&   1\\
 0&  0&   0\\
 t & 0& 0
 \end{array}
 \right)}$\ \

 \smallskip

{\rm (F1)} ${\scriptsize\left(
 \begin{array}{ccc}
 1&0&0\\
 0&0& 0\\
 0&0& 0
 \end{array}
 \right)}$\ \
{\rm (F2)} ${\scriptsize\left(
 \begin{array}{ccc}
 0 &0 &0\\
 0&0& 1\\
 0&0& 0
 \end{array}
 \right)}$\ \
{\rm (F3)} ${\scriptsize\left(
 \begin{array}{ccc}
 0&0&0\\
 0& 0& 0\\
 0&0& 1
 \end{array}
 \right)}$\ \
 {\rm (F4)} ${\scriptsize\left(
 \begin{array}{ccc}
 0 &0 &0\\
 1&0& 0\\
 0&0& 0
 \end{array}
 \right)}$\ \

 \smallskip

{\rm (F5)} ${\scriptsize\left(
 \begin{array}{ccc}
 0&0&0\\
 0& 0& 0\\
 0&0& 0
 \end{array}
 \right)}$\ \
\rr{F6} ${\scriptsize\left(
 \begin{array}{ccc}
 0&   0&   1\\
 0&  0&   0\\
 0 & 0& 0
 \end{array}
 \right)}$\ \
\end{thm}

\demo
Suppose that $G$ and $\bar{G}$ and two groups described in the theorem. By Theorem \ref{isomorphic-1} and \ref{isomorphic-2}, $\bar{G}\cong G$ if and only if there exist invertible matrices $X_2={\scriptsize\left(
 \begin{array}{ccc}
 x_{11}&      0&      0\\
 x_{21}& x_{22}&      x_{23}\\
 x_{31}& x_{32}&   x_{33}
 \end{array}
 \right)}$ and $X={\scriptsize\left(
 \begin{array}{ccc}
x_{11} & x_{12}  & x_{13}\\
0 & x_{22}  &x_{23}\\
0 &  x_{32}  & x_{33}
 \end{array}
 \right)}$ such that
$w(\bar{G})=\det(X)^{-1}X_2w(G)X^t.$

Notice that $X_2={\scriptsize\left(
 \begin{array}{ccc}
 1&      0&      0\\
 x_{21}x_{11}^{-1}& 1&      0\\
 x_{31}x_{11}^{-1} &0&1
 \end{array}
 \right)}{\scriptsize\left(
 \begin{array}{ccc}
 x_{11}&      0&      0\\
 0    & x_{22}&      x_{23}\\
 0    & x_{32}   &   x_{33}
 \end{array}
 \right)}$, and $X^t$ has same decomposition. We can use the following three kinds of operation to simplify $w(G)$.

Operation I. Take  $X_2=X={\scriptsize\left(
 \begin{array}{ccc}
 x_{11}&      0&      0\\
 0    & x_{22}&      x_{23}\\
 0    & x_{32}   &   x_{33}
 \end{array}
 \right)}$.  Then $w(\bar{G})=\det(X)^{-1}Xw(G)X^t$.

Operation II. Take $X_2={\scriptsize\left(
 \begin{array}{ccc}
 1&      0&      0\\
 x_{21}& 1&      0\\
 x_{31} & 0 &    1
 \end{array}
 \right)}$ and $X=E_3$. Then $w(\bar{G})=X_2w(G)$. This operation means adding the 1st row multiplied by a scalar $x_{21}$ to the 2nd row, and adding the 1st row multiplied by a scalar $x_{31}$ to the 3rd row.

Operation III. Take $X_2=E_3$ and $X^t={\scriptsize\left(
 \begin{array}{ccc}
 1&      0&      0\\
 x_{12}& 1&      0\\
 x_{13} &0     & 1
 \end{array}
 \right)}$. Then $w(\bar{G})=w(G)X^t$. This operation means adding the 3rd column multiplied by a scalar $x_{13}$ to the 1st column, and adding the 2nd column multiplied by a scalar $x_{12}$ to the 1st row.

According to the 1st row of $w(G)$, we have three types of $w(G)$:

  {\rm (a)} ${\scriptsize\left(
 \begin{array}{ccc}
 w_{11}& 0&0\\
 *& * &*\\
 *& *&*
 \end{array}
 \right)}$ where $w_{11}\neq 0$;\ \
 {\rm (b)} ${\scriptsize\left(
 \begin{array}{ccc}
 0& 0&0\\
 *& w_{22} &w_{23}\\
 *& w_{32}&w_{33}
 \end{array}
 \right)}$;

 {\rm (c)} ${\scriptsize\left(
 \begin{array}{ccc}
 *& w_{12}&w_{13}\\
 *& *&*\\
 *& *&*
 \end{array}
 \right)}$ where $(w_{12},w_{13})\neq (0,0)$.

Using operation I, II and III, we can not change the type of $w(G)$. Hence matrices of different type determine non-isomorphic groups.

\medskip

Case 1. $w(G)$ is of Type (a).

Using Operation I and II, $w(G)$ can be simplified to be (a$'$): $\diag(1,W)$. Obviously, two groups with characteristic matrices $w(G)$ and $w(\bar{G})$ of Type (a$'$) are mutually isomorphic if and only if there exists an invertible matrix $X=\diag(x_{11},Y)$ such that $w(\bar{G})=\det(X)^{-1}Xw(G)X^t$.
If $\diag(1,\bar{W})=\det(X)^{-1}X\diag(1,W)X^t$, then, by calculation, we have $x_{11}=\det(Y)$ and
\begin{equation}
\bar{W}=\det(Y)^{-2}YWY^t
\end{equation}
Let $Z=\det(Y)^{-1}Y$. Then $\bar{W}=ZWZ^t$. That is, $\bar{W}$ and $W$ are mutually congruent.
Using Lemma \ref{hetong-p}, \ref{hetong-2} and \ref{hetong-p-non-invertible}, we get characteristic matrices (D1)--(D7), (E8)--(E9) and (F1),
and different matrices give non-isomorphic groups.

\medskip

Case 2. $w(G)$ is of Type (b).

Using operation I, II and III, we can not change the rank of $W={\left(
 \begin{array}{cc}
 w_{22}& w_{23}\\
 w_{32}& w_{33}
 \end{array}
 \right)}$. Hence matrices with different rank of $W$ determine non-isomorphic groups.

\medskip

Subcase 2.1. $\rank(W)=2$.

Using Operation III, $w(G)$ can be simplified to be (b1): $\diag(0,W)$. Obviously, two groups with characteristic matrices $w(G)$ and $w(\bar{G})$ of Type (b1) are mutually isomorphic if and only if there exists an invertible matrix $X=\diag(x_{11},Y)$ such that $w(\bar{G})=\det(X)^{-1}Xw(G)X^t$.
If $\diag(0,\bar{W})=\det(X)^{-1}X\diag(0,W)X^t$, then, by calculation, we have
\begin{equation}
\bar{W}=x_{11}^{-1}\det(Y)^{-1}YWY^t
\end{equation}
That is, $\bar{W}$ and $W$ are mutually sub-congruent.
 Using Lemma \ref{hetong-p} and \ref{hetong-2}, we get characteristic matrices (E1)--(E7),
and different matrices give non-isomorphic groups.

\medskip

Subcase 2.2. $\rank(W)=1$.

Using Operation I, where $X={\diag(\lambda^{-1}\det(Y)^{-1},Y)}$, $w(G)$ is simplified to be ${\left(
 \begin{array}{ccc}
 0 &   0\\
 * & \lambda YWY^t
 \end{array}
 \right)}$. By Lemma \ref{hetong-p-non-invertible}, we may assume that $W=\left(
 \begin{array}{cc}
 0 & 1\\
 0 & 0
 \end{array}
 \right)$ or $\left(
 \begin{array}{cc}
 0 & 0\\
 0 & 1
 \end{array}
 \right)$. Hence there are two types:
 {\rm (b2)}: ${\scriptsize\left(
 \begin{array}{ccc}
 0& 0&0\\
 *& 0 &1\\
 *& 0&0
 \end{array}
 \right)}$ and
 {\rm (b3)} ${\scriptsize\left(
 \begin{array}{ccc}
 0& 0&0\\
 *& 0 &0\\
 *& 0&1
 \end{array}
 \right)}$.
  Using operation III, we get
(b$2'$): ${\scriptsize\left(
 \begin{array}{ccc}
 0&   0&    0\\
 0&  0&   1\\
 w_{31} & 0 & 0
 \end{array}
 \right)}$
and (b$3'$): ${\scriptsize\left(
 \begin{array}{ccc}
 0&   0&    0\\
 w_{21}&  0&   0\\
 0 & 0 & 1
 \end{array}
 \right)}$.

 Using operation I, II and III, we can not change the type of $w(G)$ from (b$2'$) to (b$3'$). Hence matrices of different type determine non-isomorphic groups.
In this subcase, we get characteristic matrices (E10)--(E11) and (F2)--(F3).

\medskip

Subcase 2.3. $W=0$.

If $(w_{21},w_{31})\neq (0,0)$, then, using operation I, we get the characteristic matrix (F4).
If $(w_{21},w_{31})=(0,0)$, then we get the characteristic matrix (F5).

\medskip

Case 3. $w(G)$ is of Type (c).

Using Operation I, we may assume that $(w_{12},w_{13})=(0,1)$. Furthermore, using Operation III, $w(G)$ can be simplified to be
Type (c$'$): ${\scriptsize\left(
 \begin{array}{ccc}
 0& 0  & 1 \\
 *& w_{22}& *\\
 * & * & *
 \end{array}
 \right)}$.
According to the value of $w_{22}$, there are two types:

(c1) ${\scriptsize\left(
 \begin{array}{ccc}
 0& 0  & 1 \\
 *& w_{22} & *\\
 * & * & *
 \end{array}
 \right)}$ where $w_{22}\neq 0$;\ \  (c2): ${\scriptsize\left(
 \begin{array}{ccc}
 0& 0  & 1 \\
 *& 0& *\\
 * & w_{32} & *
 \end{array}
 \right)}$.

We claim that two groups with characteristic matrices of Type (c1) and (c2) respectively are not isomorphic to each other. Otherwise, there exist a group whose characteristic matrix is $w(G)={\scriptsize\left(
 \begin{array}{ccc}
 0& 0  & 1 \\
 w_{21}& w_{22}& w_{23}\\
 w_{31} & w_{32} & w_{33}
 \end{array}
 \right)}$ where $w_{22}\neq 0$ ($a_1,a_2,a_3$ is a set of characteristic generators about $w(G)$), and $w(\bar{G})={\scriptsize\left(
 \begin{array}{ccc}
 0& 0  & 1 \\
 \bar{w}_{21}& 0& \bar{w}_{23}\\
 \bar{w}_{31} & \bar{w}_{32} & \bar{w}_{33}
 \end{array}
 \right)}$ is another characteristic matrix of $G$ ($\bar{a}_1,\bar{a}_2,\bar{a}_3$ is a set of characteristic generators about $w(\bar{G})$). By Theorem \ref{isomorphic-1} and \ref{isomorphic-2}, there exists invertible matrices $X_2={\scriptsize\left(
 \begin{array}{ccc}
 x_{11}&      0&      0\\
 x_{21}& x_{22}&      x_{23}\\
 x_{31}& x_{32}&   x_{33}
 \end{array}
 \right)}$ and $X={\scriptsize\left(
 \begin{array}{ccc}
x_{11} & x_{12}  & x_{13}\\
0 & x_{22}  &x_{23}\\
0 &  x_{32}  & x_{33}
 \end{array}
 \right)}$ such that
$w(\bar{G})=\det(X)^{-1}X_2w(G)X^t.$ Since the $(1,2)^{th}$ element of $X_2w(G)X^t$ is $x_{23}$, we have $x_{23}=0$ and hence $\det(X)=x_{11}x_{22}x_{33}$. By calculation, the $(2,2)^{th}$ element of $\det(X)^{-1}X_2w(G)X^t=x_{11}^{-1}x_{33}^{-1}x_{22}w_{22}\neq 0$, which is a contradiction.

\medskip

Subcase 3.1. $w(G)$ is of Type (c1).

We use following operations to simplify $w(G)$:
${\scriptsize\left(
 \begin{array}{ccc}
 0& 0  & 1 \\
 *& w_{22}& *\\
 *  & * & *
 \end{array}
 \right)}\xrightarrow{\rm {operation \ I}}{\scriptsize\left(
 \begin{array}{ccc}
 0&   0&   1\\
 *&  w_{22}&   *\\
 * & 0 & *
 \end{array}
 \right)}$
 $\xrightarrow{\rm {operation \ III}}{\scriptsize\left(
 \begin{array}{ccc}
 0& 0  & 1 \\
 0& w_{22}& *\\
 *  & 0 & *
 \end{array}
 \right)}\xrightarrow{\rm {operation \ II}}{\scriptsize\left(
 \begin{array}{ccc}
 0& 0  & 1 \\
 0& w_{22}& 0\\
 w_{31}  & 0 & 0
 \end{array}
 \right)}$
$\xrightarrow[X={\scriptsize\diag}(1,1,w_{22})]{\rm {operation \ I}}{\scriptsize\left(
 \begin{array}{ccc}
 0&   0&   1\\
 0&  1&   0\\
 w_{31} & 0 & 0
 \end{array}
 \right)}$. In this subcase, we get characteristic matrices (D8) and (E12). Similar argument as that before Subcase 3.1, we know that different $t$ for the matrix (D8) give non-isomorphic groups.

 \medskip

Subcase 3.2. $w(G)$ is of Type (c2).

If $w(G)$ is invertible, then $w_{21}\neq 0$ and $w_{32}\neq 0$. We use following operations to simplify $w(G)$:
${\scriptsize\left(
 \begin{array}{ccc}
 0& 0  & 1 \\
 w_{21}& 0& *\\
 *  & w_{32} & *
 \end{array}
 \right)}\xrightarrow{\rm {operation \ III}}{\scriptsize\left(
 \begin{array}{ccc}
 0&   0&   1\\
 w_{21}&  0&   *\\
 0 & w_{32} & *
 \end{array}
 \right)}\xrightarrow{\rm {operation \ II}}{\scriptsize\left(
 \begin{array}{ccc}
 0& 0  & 1 \\
  w_{21}&0& 0\\
 0  & w_{32} & 0
 \end{array}
 \right)}\xrightarrow[X={\scriptsize\diag}(w_{32},1,w_{21})]{\rm {operation \ I}}{\scriptsize\left(
 \begin{array}{ccc}
 0&   0&   1\\
 1&  0&   0\\
 0 & 1& 0
 \end{array}
 \right)}$. Hence we get the characteristic matrix (D9).

If $w_{21}=0$ and $w_{32}\neq 0$, then we use following operations to simplify $w(G)$:
${\scriptsize\left(
 \begin{array}{ccc}
 0& 0  & 1 \\
 0& 0& *\\
 * & w_{32} & *
 \end{array}
 \right)}$
 $\xrightarrow{\rm {operation \ III}}{\scriptsize\left(
 \begin{array}{ccc}
 0&   0&   1\\
 0&  0&   *\\
 0 & w_{32} & *
 \end{array}
 \right)}\xrightarrow{\rm {operation \ II}}{\scriptsize\left(
 \begin{array}{ccc}
 0& 0  & 1 \\
 0&0& 0\\
 0  & w_{32} & 0
 \end{array}
 \right)}\xrightarrow[X={\scriptsize\diag}(w_{32},1,1)]{\rm {operation \ I}}{\scriptsize\left(
 \begin{array}{ccc}
 0&   0&   1\\
 0&  0&   0\\
 0 & 1& 0
 \end{array}
 \right)}$. Hence we get the characteristic matrix (E13).

If $w_{21}\neq 0$ and $w_{32}=0$, then we use following operations to simplify $w(G)$:
${\scriptsize\left(
 \begin{array}{ccc}
 0& 0  & 1 \\
 w_{21}& 0& *\\
 *  & 0 & *
 \end{array}
 \right)}$
 $\xrightarrow{\rm {operation \ I}}{\scriptsize\left(
 \begin{array}{ccc}
 0&   0&   1\\
 w_{21}&  0&   *\\
 0 & 0 & *
 \end{array}
 \right)}\xrightarrow{\rm {operation \ II}}{\scriptsize\left(
 \begin{array}{ccc}
 0& 0  & 1 \\
  w_{21}&0& 0\\
 0  & 0 & 0
 \end{array}
 \right)}$
 $\xrightarrow[X={\scriptsize\diag}(1,1,w_{21})]{\rm {operation \ I}}{\scriptsize\left(
 \begin{array}{ccc}
 0&   0&   1\\
 1&  0&   0\\
 0 & 0& 0
 \end{array}
 \right)}$. Hence we get the characteristic matrix (E14).

If $w_{21}=w_{32}=0$, then, using Operation II, $w(G)$ can be simplified to be ${\scriptsize\left(
 \begin{array}{ccc}
 0&   0&   1\\
 0 &  0&   0\\
 w_{31} & 0 & 0
 \end{array}
 \right)}$. Hence we get characteristic matrices (E15) and (F6).

 Similar argument as that before Subcase 3.1, we know matrices of different types (E13)--(E15) give non-isomorphic groups, (This result also follows from Theorem \ref{1=2>3-property}, in which (H13)--(H15) have different properties.) different $t$ for the matrix (E15) give non-isomorphic groups.
 \qed

 \medskip

 We give another example about calculating $I_{\min}$ and $I_{\max}$.

 \begin{thm}\label{E4-property} Let $G$ be a finite $p$-group determined by the characteristic matrix {\rm (E4)} in Theorem $\ref{1>2=3}.$
Then $I_{\min}=m_3+1$, $I_{\max}=m_1$ for $m_1>m_2+1$ or $-r\not\in (F_p^*)^2$, and $I_{\max}=m_1+1$ for $m_1=m_2+1$ and $-r\in (F_p^*)^2$.
\end{thm}
\demo
Since $\rank(w(G))=2$, by Theorem \ref{minimal-index} (2), $m_3\le I_{\min} \le m_3+1$. By Theorem \ref{minimal-index} (1), $I_{\min} =m_3$ if and only if there exist invertible matrices $X_2={\scriptsize\left(
 \begin{array}{ccc}
 x_{11}&      0&      0\\
 x_{21}& x_{22}&      x_{23}\\
 x_{31}& x_{32}&   x_{33}
 \end{array}
 \right)}$ and $X={\scriptsize\left(
 \begin{array}{ccc}
x_{11} & x_{12}  & x_{13}\\
0 & x_{22}  &x_{23}\\
0 &  x_{32}  & x_{33}
 \end{array}
 \right)}$ and
$w(\bar{G})=\det(X)^{-1}X_2w(G)X^t=(\bar{w}_{ij})$ such that $\left(
 \begin{array}{cc}
 \bar{w}_{11}& \bar{w}_{12}\\
 \bar{w}_{21} &\bar{w}_{22}
 \end{array}
 \right)$ is invertible. Let $Y= \left(
 \begin{array}{cc}
 x_{22}& x_{23}\\
 x_{32} &x_{33}
 \end{array}
 \right)$ and $W=\left(
 \begin{array}{cc}
 1& 1\\
 -1 &r
 \end{array}
 \right)$. By calculation, we have
 \begin{eqnarray*}\label{E4-1}
 \det(X)w(\bar{G})&=&\left(
 \begin{array}{cc}
 0&      O\\
 YW(x_{12},x_{13})^t& YWY^t
 \end{array}
 \right)\\
 &=&{\scriptsize\left(
 \begin{array}{ccc}
0 & 0  & 0\\
* & *  &*\\
x_{12}x_{32}-x_{12}x_{33}+rx_{13}x_{33}+x_{13}x_{32} &  *  & x_{32}^2+rx_{33}^2
 \end{array}
 \right)}
 \end{eqnarray*}

 Since $\left(
 \begin{array}{cc}
 \bar{w}_{11}& \bar{w}_{12}\\
 \bar{w}_{21} &\bar{w}_{22}
 \end{array}
 \right)$ is not invertible, we have $I_{\min}=m_3+1$.

Since $YWY^t$ is invertible, by Theorem \ref{maximal-index}, $I_{\max}\neq m_1+2$ and Condition (a) in Theorem \ref{maximal-index} (4) does not hold. If $-r$ is not a square, then $x_{32}+rx_{33}\neq 0$. Hence Condition (b) in Theorem \ref{maximal-index} (4) does not hold. In this case, $I_{\max}=m_1$. If $-r=s^2$ is a square, then, taking $x_{12}=x_{32}=sx_{33}=sx_{13}$, we have $\left(
 \begin{array}{cc}
 \bar{w}_{11}& \bar{w}_{13}\\
 \bar{w}_{31} &\bar{w}_{33}
 \end{array}
 \right)=0$. Hence Condition (b) in Theorem \ref{maximal-index} (4) holds if and only if $m_1=m+2+1$. In this case, $I_{\max}=m_1$ for $m_1>m_2+1$, and $I_{\max}=m_1+1$ for $m_1=m_2+1$.\qed

\begin{thm}\label{1>2=3-property} Let $G$ be a finite $p$-group determined by a characteristic matrix in Theorem $\ref{1>2=3}.$

\rr{1} If $I_{\min}={m_3}$, then the characteristic matrix is one of the following: {\rm (D2)--(D4), (D6)--(D9), (E8)--(E9), (E12)--(E15)};

\rr{2} If $I_{\min}={m_3+1}$, then the characteristic matrix is one of the following: {\rm (D1), (D5), (E1)--(E7), (E10)--(E11), (F1)--(F4), (F6)};

\rr{3} If $I_{\min}={m_3+2}$, then the characteristic matrix is {\rm (F5)};

\rr{4} If $I_{\max}={m_1+2}$, then the characteristic matrix is one of the following: {\rm (E14)--(E15), (F1), (F4)--(F6)};

\rr{5} If $I_{\max}={m_1+1}$, then the characteristic matrix is one of the following: {\rm (D8)--(D9), (E1)--(E2) for $m_1=m_2+1$, (E3) for $m_1=m_2+1$ and $-\nu\in (F_p^*)^2$, (E4) for $m_1=m_2+1$ and $-r\in (F_p^*)^2$, (E5)--(E6) for $m_1=m_2+1$, (E8)--(E13), (F2)--(F3)};

\rr{6} If $I_{\max}={m_1}$, then the characteristic matrix is one of the following: {\rm (D1)--(D7), (E1)--(E2) for $m_1>m_2+1$, (E3) for $m_1>m_2+1$ or $-\nu\not\in (F_p^*)^2$, (E4) for $m_1>m_2+1$ or $-r\not\in (F_p^*)^2$, (E5)--(E6) for $m_1>m_2+1$, (E7)}.

\end{thm}

 \section{The case of $m_1=m_2>m_3$}

In this section, $p$ is either odd or $2$.

\begin{thm}\label{1=2>3} Let $G$ be a finite $p$-group with $d(G)=3$, $\Phi(G)\le Z(G)$ and $G'\cong C_p^3$. If the type of $G/G'$ is $(p^{m_1},p^{m_2},p^{m_3})$ where $m_1=m_2>m_3$, then we can choose suitable generators of $G$, such that the characteristic matrix of $G$ is one of the
following matrices, and different matrices give non-isomorphic groups. Where $\eta$ is a fixed square
non-residue modulo odd $p$, $\nu=1$ or $\eta$, $t\neq 0$, and  $r=1,2,\dots,p-2$.

\smallskip

\rr{1} For odd $p$

\medskip

\rr{G1} ${\scriptsize\left(
 \begin{array}{ccc}
 0 & -1 &0\\
 1&  0& 0\\
 0& 0& 1
 \end{array}
 \right)}$\ \
{\rm (G2)} ${\scriptsize\left(
 \begin{array}{ccc}
\nu&-1&0\\
1& 0& 0\\
 0& 0& 1
 \end{array}
 \right)}$\ \
{\rm (G3)} ${\scriptsize\left(
 \begin{array}{ccc}
 1  & 0  & 0\\
 0 &\nu& 0\\
 0 &0 & 1
 \end{array}
 \right)}$\ \
{\rm (G4)} ${\scriptsize\left(
 \begin{array}{ccc}
 1  & -1  & 0\\
 1 &r& 0\\
 0 &0& 1
 \end{array}
 \right)}$\ \

\smallskip

\rr{H1} ${\scriptsize\left(
 \begin{array}{ccc}
 0 & -1 &0\\
 1&  0& 0\\
 0& 0& 0
 \end{array}
 \right)}$\ \
{\rm (H2)} ${\scriptsize\left(
 \begin{array}{ccc}
1&-1&0\\
1& 0& 0\\
 0& 0& 0
 \end{array}
 \right)}$\ \
{\rm (H3)} ${\scriptsize\left(
 \begin{array}{ccc}
 1  & 0  & 0\\
 0 &\nu & 0\\
 0 &0 & 0
 \end{array}
 \right)}$\ \
{\rm (H4)} ${\scriptsize\left(
 \begin{array}{ccc}
 1 &-1& 0 \\
 1 &r& 0\\
 0 &0& 0
 \end{array}
 \right)}$

\medskip

\rr{2} For $p=2$

\medskip

{\rm (G5)} ${\scriptsize\left(
 \begin{array}{ccc}
0&1&0\\
1& 0& 0\\
 0& 0& 1
 \end{array}
 \right)}$\ \
\rr{G6} ${\scriptsize\left(
 \begin{array}{ccc}
 1 & 0 &0\\
 0&  1& 0\\
 0& 0& 1
 \end{array}
 \right)}$\ \
{\rm (G7)} ${\scriptsize\left(
 \begin{array}{ccc}
 1  & 1  & 0\\
 0 &1& 0\\
 0 &0 & 1
 \end{array}
 \right)}$

\smallskip

{\rm (H5)} ${\scriptsize\left(
 \begin{array}{ccc}
0&1&0\\
1& 0& 0\\
 0& 0& 0
 \end{array}
 \right)}$\ \
 \rr{H6} ${\scriptsize\left(
 \begin{array}{ccc}
 1 & 0 &0\\
 0&  1& 0\\
 0& 0& 0
 \end{array}
 \right)}$\ \
{\rm (H7)} ${\scriptsize\left(
 \begin{array}{ccc}
 1  & 1  & 0\\
 0 &1& 0\\
 0 &0 & 0
 \end{array}
 \right)}$

 \medskip

 \rr{3} For both odd $p$ and $p=2$

 \medskip

\rr{G8} ${\scriptsize\left(
 \begin{array}{ccc}
 1&   0&   0\\
 0&  0&   1\\
 0 & t & 0
 \end{array}
 \right)}$\ \
 \rr{G9} ${\scriptsize\left(
 \begin{array}{ccc}
 0&   1&   0\\
 0&  0&   1\\
 1 & 0& 0
 \end{array}
 \right)}$\ \

\smallskip

{\rm (H8)} ${\scriptsize\left(
 \begin{array}{ccc}
 0 &0 &0\\
 1&0& 0\\
 0&0& 1
 \end{array}
 \right)}$\ \
{\rm (H9)} ${\scriptsize\left(
 \begin{array}{ccc}
 0&0&0\\
 0& \nu& 0\\
 0&0& 1
 \end{array}
 \right)}$\ \
 \rr{H10} ${\scriptsize\left(
 \begin{array}{ccc}
 0&   0&    0\\
 1&  0&   0\\
 0 & 1 & 0
 \end{array}
 \right)}$\ \
\rr{H11} ${\scriptsize\left(
 \begin{array}{ccc}
 0&   0&    0\\
 0&  1&   0\\
 1 & 0 & 0
 \end{array}
 \right)}$\ \

 \smallskip

\rr{H12} ${\scriptsize\left(
 \begin{array}{ccc}
 1&   0&   0\\
 0&  0&   1\\
 0 & 0 & 0
 \end{array}
 \right)}$\ \
 \rr{H13} ${\scriptsize\left(
 \begin{array}{ccc}
 0&   1&   0\\
 0&  0&   1\\
 0 & 0& 0
 \end{array}
 \right)}$\ \
 \rr{H14} ${\scriptsize\left(
 \begin{array}{ccc}
 0&   0&   0\\
 0&  0&   1\\
 1 & 0& 0
 \end{array}
 \right)}$\ \
\rr{H15} ${\scriptsize\left(
 \begin{array}{ccc}
 0&   0&   0\\
 0&  0&   1\\
 0 & t& 0
 \end{array}
 \right)}$\ \

 \smallskip

{\rm (I1)} ${\scriptsize\left(
 \begin{array}{ccc}
 0&0&0\\
 0&0& 0\\
 0&0& 1
 \end{array}
 \right)}$\ \
{\rm (I2)} ${\scriptsize\left(
 \begin{array}{ccc}
 0 &0 &0\\
 1&0& 0\\
 0&0& 0
 \end{array}
 \right)}$\ \
{\rm (I3)} ${\scriptsize\left(
 \begin{array}{ccc}
 0&0&0\\
 0&1& 0\\
 0&0& 0
 \end{array}
 \right)}$\ \
 {\rm (I4)} ${\scriptsize\left(
 \begin{array}{ccc}
 0 &0 &0\\
 0&0& 0\\
 1&0& 0
 \end{array}
 \right)}$\ \

 \smallskip

{\rm (I5)} ${\scriptsize\left(
 \begin{array}{ccc}
 0&0&0\\
 0& 0& 0\\
 0&0& 0
 \end{array}
 \right)}$\ \
\rr{I6} ${\scriptsize\left(
 \begin{array}{ccc}
 0&   0&   0\\
 0&  0&   1\\
 0 & 0& 0
 \end{array}
 \right)}$\ \
\end{thm}
\demo  Suppose that $G$ and $\bar{G}$ and two groups described in the theorem. By Theorem \ref{isomorphic-1} and \ref{isomorphic-2}, $\bar{G}\cong G$ if and only if there exists invertible matrices $X_2={\scriptsize\left(
 \begin{array}{ccc}
 x_{11}&      x_{12}&      0\\
 x_{21}& x_{22}&      0\\
 x_{31}& x_{32}&   x_{33}
 \end{array}
 \right)}$ and $X={\scriptsize\left(
 \begin{array}{ccc}
x_{11} & x_{12}  & x_{13}\\
x_{21} & x_{22}  &x_{23}\\
0 &  0  & x_{33}
 \end{array}
 \right)}$ such that
$w(\bar{G})=\det(X)^{-1}X_2w(G)X^t.$

Let $P={\scriptsize\left(
 \begin{array}{ccc}
0 & 0  & 1\\
1 & 0  & 0\\
0 & 1  & 0
 \end{array}
 \right)}$. Then $P^{-1}=P^t={\scriptsize\left(
 \begin{array}{ccc}
0 & 1  & 0\\
0 & 0  & 1\\
1 & 0  & 0
 \end{array}
 \right)}$. Let $X_2'=PXP^{-1}$ and $X'=PX_2P^{-1}$. Then $Pw(\bar{G})^tP^{-1}=\det(X')^{-1}X_2'(Pw(G)^tP^{-1})(X')^t$. Since $X_2'={\scriptsize\left(
 \begin{array}{ccc}
 x_{33}&      0&      0\\
 x_{13}& x_{11}&      x_{12}\\
 x_{23}& x_{21}&   x_{22}
 \end{array}
 \right)}$ and $X'={\scriptsize\left(
 \begin{array}{ccc}
x_{33} & x_{31}  & x_{32}\\
0 & x_{11}  & x_{12}\\
0 & x_{21}  & x_{22}
 \end{array}
 \right)}$, $Pw(\bar{G})^tP^{-1}$ and $Pw(G)^tP^{-1}$ have relation in Theorem \ref{1>2=3}. So, compared with last section, we get the result listed in the theorem.\qed

\begin{thm}\label{1=2>3-property} Let $G$ be a finite $p$-group determined by a characteristic matrix in Theorem $\ref{1=2>3}.$

\rr{1} If $I_{\min}={m_3}$, then the characteristic matrix is one of the following: {\rm (G1)--(G9), (H1)--(H7), (H12)--(H13)};

\rr{2} If $I_{\min}={m_3+1}$, then the characteristic matrix is one of the following: {\rm (H8)--(H11), (H14)--(H15), (I2)--(I3), (I6)};

\rr{3} If $I_{\min}={m_3+2}$, then the characteristic matrix is one of the following: {\rm (I1), (I4)--(I5)};

\rr{4} If $I_{\max}={m_1+2}$, then the characteristic matrix is one of the following: {\rm (H1)--(H2),(H3) for $-\nu\in (F_p^*)^2$, (H4) for $-r\in (F_p^*)^2$, (H5)--(H6), (H10), (H13), (H15), (I2)--(I6)};

\rr{5} If $I_{\max}={m_1+1}$, then the characteristic matrix is one of the following: {\rm (G1)--(G2), (G3) for $-\nu\in (F_p^*)^2$, (G4) for $-r\in (F_p^*)^2$, (G5)--(G6), (G8)--(G9), (H3) for $-\nu\not\in (F_p^*)^2$, (H4) for $-r\not\in (F_p^*)^2$, (H7)--(H9), (H11)--(H12), (H14), (I1)};

\rr{6} If $I_{\max}={m_1}$, then the characteristic matrix is one of the following: {\rm (G3) for $-\nu\not\in (F_p^*)^2$, (G4) for $-r\not\in (F_p^*)^2$, (G7)}.

\end{thm}

 \section{The case of $m_1=m_2=m_3$ for odd prime $p$.}

\begin{thm}\label{1=2=3}
 Let $G$ be a finite $p$-group with $d(G)=3$, $\Phi(G)\le Z(G)$ and $G'\cong C_p^3$, where $p$ is an odd prime. If the type of $G/G'$ is $(p^{m_1},p^{m_2},p^{m_3})$ where $m_1=m_2=m_3$, then we can choose suitable generators of $G$, such that the characteristic matrix of $G$ is one of the
following matrices, where different types of matrix give non-isomorphic groups. $($Where $\eta$ is a fixed square
non-residue modulo $p$, $\nu=1$ or $\eta$, and $r=1,2,\dots,p-2. )$

\medskip

\noindent{\rm (J1)} ${\scriptsize\left(
 \begin{array}{ccc}
 1& 0&0\\
 0& 1&0\\
 0& 0&1
 \end{array}
 \right)}$\ \
{\rm (J2)} ${\scriptsize\left(
 \begin{array}{ccc}
 1& 0&0\\
 0& 0&1\\
 0& -1&0
 \end{array}
 \right)}$\ \
  {\rm (J3)} ${\scriptsize\left(
 \begin{array}{ccc}
 1& 0&0\\
 0& \nu & 1\\
 0& -1& 0
 \end{array}
 \right)}$\ \
{\rm (J4)} ${\scriptsize\left(
 \begin{array}{ccc}
 1& 0&0\\
 0& 1&1\\
 0& -1&r
 \end{array}
 \right)}$\ \
{\rm (J5)} ${\scriptsize\left(
 \begin{array}{ccc}
 0& 0&1\\
 0& 1&1\\
 1& -1&0
 \end{array}
 \right)}$\ \

\smallskip

\noindent {\rm (K1)} ${\scriptsize\left(
 \begin{array}{ccc}
 1& 0&0\\
 0& \nu &0\\
 0& 0&0
 \end{array}
 \right)}$\ \
{\rm (K2)} ${\scriptsize\left(
 \begin{array}{ccc}
 1& 0&0\\
 0& 0& 1\\
 0& 0& 0
 \end{array}
 \right)}$\ \
 {\rm (K3)} ${\scriptsize\left(
 \begin{array}{ccc}
 0& 0&1\\
 0& 0&0\\
 0& 1& 0
 \end{array}
 \right)}$\ \
  {\rm (K4)} ${\scriptsize\left(
 \begin{array}{ccc}
 0& 0&0\\
 0& 0&1\\
 0& -1&0
 \end{array}
 \right)}$\ \
{\rm (K5)} ${\scriptsize\left(
 \begin{array}{ccc}
 0& 0&0\\
 0& 1&1\\
 0& -1&0
 \end{array}
 \right)}$\ \

\smallskip

\noindent {\rm (K6)} ${\scriptsize\left(
 \begin{array}{ccc}
 0& 0&0\\
 0& 1&1\\
 0& -1&r
 \end{array}
 \right)}$\ \ \
\noindent{\rm (L1)} ${\scriptsize\left(
 \begin{array}{ccc}
 0& 0&0\\
 0& 0&0\\
 0& 0&0
 \end{array}
 \right)}$\ \
 {\rm (L2)} ${\scriptsize\left(
 \begin{array}{ccc}
 1& 0&0\\
 0& 0& 0\\
 0& 0& 0
 \end{array}
 \right)}$\ \
 {\rm (L3)} ${\scriptsize\left(
 \begin{array}{ccc}
 0& 0&0\\
 0& 0& 1\\
 0& 0& 0
 \end{array}
 \right)}$\ \

 \end{thm}
\demo Suppose that $G$ and $\bar{G}$ and two groups described in the theorem. By Theorem \ref{isomorphic-1}, $\bar{G}\cong G$ if and only if there exists invertible matrix $X={\scriptsize\left(
 \begin{array}{ccc}
 x_{11}& x_{12}&      x_{13}\\
 x_{21}& x_{22}&      x_{23}\\
 x_{31}& x_{32}&      x_{33}
 \end{array}
 \right)}$ such that
$w(\bar{G})=\det(X)^{-1}Xw(G)X^t.$ For short, we call $w(\bar{G})$ and $w(G)$ are quasi-congruent.

Case 1: $w(G)$ is symmetric.

By an
elementary property of symmetric matrices over fields of
characteristic different from 2, $w(G)$ is congruent to a diagonal matrix. Hence we may let
$w(G)=\diag(i,s,w)$. If $w(G)=0$, then we get the matrix (L1). If $w(G)\neq 0$, then, without loss of generality, we may assume that $i\neq 0$.
Take $X=\diag(i^{-1}\det(Y)^{-1},\det(Y)^{-1}Y)$. Then $w(\bar{G})=\det(X)^{-1}Xw(G)X^t=\diag(1,Y\diag(is,iw)Y^t)$.
By Lemma \ref{hetong-p} and \ref{hetong-p-non-invertible}, $w(\bar{G})$ is one of the matrices (J1), (K1), (L2) or $\diag(1,1,\eta)$.

If $w(\bar{G})=\diag(1,1,\eta)$, then we can simplify it to be the matrix (J1). Assume that $\eta=z^2\zeta$, where $\zeta\in F_p^*$ such that
$\zeta\not\in (F_p^*)^2$ and $\zeta-1=\gamma^2\in (F_p^*)^2$. Then, taking $X={\scriptsize\left(
 \begin{array}{ccc}
 0& 0&-1\\
 z& -z\gamma& 0\\
 -z\gamma & -z& 0
 \end{array}
 \right)}$, we have $\det(X)=\eta$ and $\det(X)^{-1}Xw(\bar{G})X^t=\diag(1,1,1)$. Hence we also get the matrix (J1).

 Obviously matrices of different types (J1), (K1), (L1) and (L2) are not quasi-congruent. We claim that different $\nu$ for the matrix (K1) give matrices which are not quasi-congruent. Otherwise, there exists invertible matrix $X={\scriptsize\left(
 \begin{array}{ccc}
x_{11} & x_{12}  & x_{13}\\
x_{21} & x_{22}  &x_{23}\\
x_{31} &  x_{32}  & x_{33}
 \end{array}
 \right)}$ such that
\begin{equation*}
\diag(1,\eta,0)=\det(X)^{-1}X\diag(1,1,0)X^t.
\end{equation*} Let $Y={\left(
 \begin{array}{cc}
x_{11} & x_{12}  \\
x_{21} & x_{22}  \\
 \end{array}
 \right)}$. By calculation, we have $\diag(1,\eta)=\det(X)^{-1}YY^t$, which contradicts Lemma \ref{hetong-p}.

\medskip

Case 2: $w(G)$ is not symmetric.

 For a characteristic matrix $w(G)$, let
$W_1=2^{-1}(w(G)+w(G)^t)$  and  $W_2=2^{-1}(w(G)-w(G)^t).$ Then $W_1$ is a
symmetric matrix and $W_2$ a skew-symmetric matrix.
Assume that $W_2=
{\scriptsize\left(
 \begin{array}{ccc}
 0    & x          &y\\
 -x   & 0          &z\\
 -y   &-z          &0
 \end{array}
 \right)}$. Let $X$ be
 ${\scriptsize\left(
 \begin{array}{ccc}
 z    & -y          &x\\
 y^{-1}    & 0                 &0\\
 0    &      0                 &1
 \end{array}
 \right)}\  {\rm for}\ y\neq 0,$
$ {\scriptsize \left(
 \begin{array}{ccc}
 z    & -y          &x\\
 0    & z^{-1}     &0\\
 0    & 0           &1
 \end{array}
 \right)}\ {\rm for}\ z\neq 0$,
 and
${\scriptsize\left(
 \begin{array}{ccc}
 z    & -y            &x\\
 x^{-1}   & 0             &1\\
 0    & 0                 &0
 \end{array}
 \right)} \ {\rm for}\ x\neq 0$ respectively. Then $\det(X)=1$ and
$\det(X)^{-1}XW_2X^t={\scriptsize\left(
 \begin{array}{ccc}
 0    & 0   &0\\
 0    & 0   &1\\
 0    &-1   &0
 \end{array}
 \right)}.$ Hence we may assume that $W_2={\scriptsize\left(
 \begin{array}{ccc}
 0    & 0   &0\\
 0    & 0   &1\\
 0    &-1   &0
 \end{array}
 \right)}$.
 Let $W_1={\scriptsize\left(
 \begin{array}{ccc}
 i    & j   &k\\
 j    & s   &t\\
 k    & t   &w
 \end{array}
 \right)}.$ Then $w(G)={\scriptsize\left(
 \begin{array}{ccc}
 i    & j   &k\\
 j    & s   &t+1\\
 k    & t-1   &w
 \end{array}
 \right)}$.
In the following, when we use $X$ to simplify $w(G)$, we hope that $X$ satisfy $\det(X)^{-1}XW_2X^t=W_2$.
By calculation, this means that $x_{11}=1$ and $x_{12}=x_{13}=0$. That is, $X={\scriptsize\left(
 \begin{array}{ccc}
1 & 0  & 0\\
x_{21} & x_{22}  &x_{23}\\
x_{31} &  x_{32}  & x_{33}
 \end{array}
 \right)}$. According to the value of $i$, there are two types of $w(G)$:

 \rr{a} ${\scriptsize\left(
 \begin{array}{ccc}
 i    & j   &k\\
 j    & s   &t+1\\
 k    & t-1   &w
 \end{array}
 \right)}$ where $i\neq 0$;\ \
 \rr{b} ${\scriptsize\left(
 \begin{array}{ccc}
 0    & j   &k\\
 j    & s   &t+1\\
 k    & t-1   &w
 \end{array}
 \right)}$

Notice that the $(1,1)^{th}$ element of $Xw(G)X^t$ is also $i$. Matrices of different types are not quasi-congruent.

\medskip

Subcase 2.1: $w(G)$ is of Type (a).

Let $X={\scriptsize\left(
 \begin{array}{ccc}
 1    & 0   &0\\
 -j    & i   &0\\
 -i^{-1}k    & 0   &1
 \end{array}
 \right)}$. Then $\det(X)^{-1}Xw(G)X^t=\diag(1,W)$. By calculation, matrices $\diag(1,W)$ and $\diag(1,\bar{W})$ are quasi-congruent if and only if there exists $Y=\left(
 \begin{array}{cc}
 x_{22}   &x_{23}\\
 x_{32}   &x_{33}\\
 \end{array}
 \right)$ such that $\bar{W}=\det(Y)^{-2}YWY^t$. Let $Z=\det(Y)^{-1}Y$. Then we have $\bar{W}=ZWZ^t$. That is, $\bar{W}$ and $W$ are mutually congruent. Using Lemma \ref{hetong-p} and \ref{hetong-p-non-invertible}, we get characteristic matrices (J2)--(J4) and (K2), and different matrices give non-isomorphic groups.

\medskip

Subcase 2.2: $w(G)$ is of Type (b).

According to the value of $(j,k)$, there are two types of $w(G)$:

 \rr{b1} ${\scriptsize\left(
 \begin{array}{ccc}
 0    & j   &k\\
 j    & s   &t+1\\
 k    & t-1   &w
 \end{array}
 \right)}$ where $(j,k)\neq (0,0)$;\ \
  \rr{b2} ${\scriptsize\left(
 \begin{array}{ccc}
 0    & 0   &0\\
 0    & s   &t+1\\
 0    & t-1   &w
 \end{array}
 \right)}$.

It is easy to see (by calculation) that matrices of different types are not quasi-congruent.

\medskip

Subcase 2.2.1: $w(G)$ is of Type (b1).

Subcase 2.2.1.1: $j=0$ and $k\neq 0$.

If $s\neq 0$, then, letting $X={\scriptsize\left(
 \begin{array}{ccc}
 1    & 0   &0\\
 -t    & k   &0\\
 -2^{-1}sw    & 0   &sk
 \end{array}
 \right)}$, $w(\bar{G})=\det(X)^{-1}Xw(G)X^t={\scriptsize\left(
 \begin{array}{ccc}
 0    & 0   &1\\
 0    & 1   &1\\
 1    & -1   &0
 \end{array}
 \right)}$. Hence we get the characteristic matrix (J5).

If $s=0$, then, letting $X={\scriptsize\left(
 \begin{array}{ccc}
 1    & 0   &0\\
 -t    & k   &0\\
 -2^{-1}k^{-1}w    & 0   &1
 \end{array}
 \right)}$, $w(\bar{G})=\det(X)^{-1}Xw(G)X^t={\scriptsize\left(
 \begin{array}{ccc}
 0    & 0   &1\\
 0    & 0   &1\\
 1    & -1   &0
 \end{array}
 \right)}$. Again let $X={\scriptsize\left(
 \begin{array}{ccc}
 -1    & -1   &0\\
 -1    & 1   &0\\
 0    & 0   &1
 \end{array}
 \right)}$, $w(\bar{G})$ can be simplified to be matrix (K3).

\medskip

Subcase 2.2.1.2: $j\neq 0$ and $k\neq 0$.

Let $X={\scriptsize\left(
 \begin{array}{ccc}
 1    & 0   &0\\
 0    & 1   &-k^{-1}j\\
 0    & 0   &1
 \end{array}
 \right)}$. Then $w(\bar{G})=\det(X)^{-1}Xw(G)X^t={\scriptsize\left(
 \begin{array}{ccc}
 0    & 0   &k\\
 0    & s'   &t'+1\\
 k    & t'-1   &w
 \end{array}
 \right)}$, where $t'=t-k^{-1}jw$. It is reduced to subcase 2.2.1.1.

\medskip

Subcase 2.2.1.3: $k=0$ and $j\neq 0$.

Let $X={\scriptsize\left(
 \begin{array}{ccc}
 1    & 0   &0\\
 0    & 1   &0\\
 0    & 1   &1
 \end{array}
 \right)}$. Then $w(\bar{G})=\det(X)^{-1}Xw(G)X^t={\scriptsize\left(
 \begin{array}{ccc}
 0    & j   &j\\
 j    & s   &t'+1\\
 j    & t'-1   &w'
 \end{array}
 \right)}$, where $t'=t+s$ and $w'=w+2t+s$. It is reduced to subcase 2.2.1.2.

\medskip

Subcase 2.2.2: $w(G)$ is of Type (b2).

In this subcase, $w(G)=\diag(0,W)$. By calculation, matrices $\diag(0,W)$ and $\diag(0,\bar{W})$ are quasi-congruent if and only if there exist $x_{11}\neq 0$ and $Y=\left(
 \begin{array}{cc}
 x_{22}   &x_{23}\\
 x_{32}   &x_{33}\\
 \end{array}
 \right)$ such that $\bar{W}=x_{11}^{-1}\det(Y)^{-1}YWY^t$. That is, $\bar{W}$ and $W$ sub-congruent. By lemma \ref{hetong-p} and \ref{hetong-p-non-invertible}, we get characteristic matrices (K4)--(K6) and (L3), and different matrices give non-isomorphic groups.
\qed

 \begin{thm}\label{J4-property} Let $G$ be a finite $p$-group determined by the characteristic matrix {\rm (J4)} in Theorem $\ref{1=2=3}.$
Then $I_{\min}=m_3$ and $I_{\max}=m_1+1$.
\end{thm}
\demo Since $\left(
 \begin{array}{cc}
 {w}_{11}& {w}_{12}\\
 {w}_{21} &{w}_{22}
 \end{array}
 \right)$ is invertible,
 by Theorem \ref{minimal-index} (1), $I_{\min} =m_3$.
 Since $\rank(w(G))=3$, by Theorem \ref{maximal-index} (3), $m_1\le I_{\max}\le m_1+1$. We will show that Condition (a) in Theorem \ref{maximal-index} (4) holds. Hence $I_{\max}= m_1+1$.

If $-r\in (F_p^*)^2$, then,
 letting $-r=y^2$ and $X={\scriptsize\left(
 \begin{array}{ccc}
 0&      1&      0\\
 1& 0&      0\\
 0& y&   1
 \end{array}
 \right)}$,
$w(\bar{G})=\det(X)^{-1}Xw(G)X^t={\scriptsize\left(
 \begin{array}{ccc}
-1 & 0  & -y-1\\
0 & -1  &0\\
1-y &  0  & 0
 \end{array}
 \right)}$. Hence Condition (a) in Theorem \ref{maximal-index} (4) holds.

If $-r\not\in (F_p^*)^2$, then, since $F_p^2\cap (-r-F_p^2)\neq \phi$, there exist $x,y\in F_p^*$ such that $x^2=-r-y^2$.
Let $X={\scriptsize\left(
 \begin{array}{ccc}
 1&      0&      0\\
  x^{-1}(1-y)& 1&      0\\
 x& y&   1
 \end{array}
 \right)}$. Then
$w(\bar{G})=\det(X)^{-1}Xw(G)X^t={\scriptsize\left(
 \begin{array}{ccc}
1 &  x^{-1}(1-y)  & x\\
 x^{-1}(1-y) &  x^{-2}(1-y)^2+1  &2\\
x &  0  & 0
 \end{array}
 \right)}$. Hence Condition (a) in Theorem \ref{maximal-index} (4) holds.
\qed

\begin{thm}\label{1=2>3-property} Let $G$ be a finite $p$-group determined by a characteristic matrix in Theorem $\ref{1=2>3}.$

\rr{1} If $I_{\min}={m_3}$, then the characteristic matrix is one of the following: {\rm (J1)--(J5), (K1)--(K6)};

\rr{2} If $I_{\min}={m_3+1}$, then the characteristic matrix is one of the following: {\rm (L2)--(L3)};

\rr{3} If $I_{\min}={m_3+2}$, then the characteristic matrix is {\rm (L1)};

\rr{4} If $I_{\max}={m_1+2}$, then the characteristic matrix is one of the following: {\rm (K1) for $-\nu\in (F_p^*)^2$, (K3)--(K5), (K6) for $-r\in (F_p^*)^2$, (L1)--(L3)};

\rr{5} If $I_{\max}={m_1+1}$, then the characteristic matrix is one of the following: {\rm (J1)--(J5), (K1) for $-\nu\in (F_p^*)^2$, (K2), (K6) for $-r\in (F_p^*)^2$}.
\end{thm}

\section{Other cases for $p=2$}

There are still three cases for $p=2$. That is, (a) $m_1=m_2=m_3=1$, (b) $m_1>m_2=m_3=1$ and (c) $m_1=m_2=m_3\ge 2$. If $m_1=m_2=m_3=1$, then $|G|=2^6$ and we can pick up what we need from the list of groups of order $2^6$. If $m_1>m_2=m_3=1$ or $m_1=m_2=m_3\ge 2$, then we need find characteristic matrices for all non-isomorphic groups. On the other hand, if we know all non-isomorphic groups, then we can get their characteristic matrices. In (b) and (c), when we consider the minimal possible value of $m_i$, we get $|G|=2^7$ for (b) and $|G|=2^9$ for (c). Fortunately, 2-groups of order less than $2^{10}$ have been classified, which can be listed by using Magma. Hence we get the following theorem \ref{1>2=3 for m_2=1 and p=2} and \ref{1=2=3 for p=2}.

\begin{thm}\label{1>2=3 for m_2=1 and p=2} Let $G$ be a finite $2$-group with $d(G)=3$, $\Phi(G)\le Z(G)$ and $G'\cong C_2^3$. If the type of $G/G'$ is $(2^{m_1},2,2)$ where $m_1>1$, then we can choose suitable generators of $G$, such that the characteristic matrix of $G$ is one of the
following matrices, and different matrices give non-isomorphic groups.

\smallskip

{\rm (M1)} ${\scriptsize\left(
 \begin{array}{ccc}
1&0&0\\
0& 0& 1\\
 0& 1& 0
 \end{array}
 \right)}$\ \
\rr{M2} ${\scriptsize\left(
 \begin{array}{ccc}
 1 & 0 &0\\
 0&  1& 0\\
 0& 0& 1
 \end{array}
 \right)}$\ \
{\rm (M3)} ${\scriptsize\left(
 \begin{array}{ccc}
 1  & 0  & 0\\
 0 &1& 0\\
 0 &1 & 1
 \end{array}
 \right)}$\ \
\rr{M4} ${\scriptsize\left(
 \begin{array}{ccc}
 0&   0&   1\\
 0&  1&   0\\
 1 & 0 & 0
 \end{array}
 \right)}$\ \

\smallskip

 \rr{M5} ${\scriptsize\left(
 \begin{array}{ccc}
 0&   0&   1\\
 1&  0&   0\\
 0 & 1& 0
 \end{array}
 \right)}$\ \
 \rr{M6} ${\scriptsize\left(
 \begin{array}{ccc}
 0&   0&   1\\
 0&  1&   0\\
 1 & 1& 0
 \end{array}
 \right)}$\ \
 \rr{M7} ${\scriptsize\left(
 \begin{array}{ccc}
 0&   1&   1\\
 0&  1&   0\\
 1 & 1& 0
 \end{array}
 \right)}$

\smallskip

{\rm (N1)} ${\scriptsize\left(
 \begin{array}{ccc}
0&0&0\\
0& 0& 1\\
 0& 1& 0
 \end{array}
 \right)}$\ \
 \rr{N2} ${\scriptsize\left(
 \begin{array}{ccc}
 0 & 0 &0\\
 0&  1& 0\\
 0& 0& 1
 \end{array}
 \right)}$\ \
{\rm (N3)} ${\scriptsize\left(
 \begin{array}{ccc}
 0  & 0  & 0\\
 0 &1& 0\\
 0 &1 & 1
 \end{array}
 \right)}$\ \
 \rr{N4} ${\scriptsize\left(
 \begin{array}{ccc}
 0&   0&   0\\
 0&  1&   1\\
 1 & 1& 1
 \end{array}
 \right)}$

\smallskip

{\rm (N5)} ${\scriptsize\left(
 \begin{array}{ccc}
 1 &0 &0\\
 0&0& 1\\
 0&0& 0
 \end{array}
 \right)}$\ \
{\rm (N6)} ${\scriptsize\left(
 \begin{array}{ccc}
 1&0&0\\
 0& 0& 0\\
 0&0& 1
 \end{array}
 \right)}$\ \
 \rr{N7} ${\scriptsize\left(
 \begin{array}{ccc}
 0&   0&    0\\
 0&  0&   1\\
 1 & 0 & 0
 \end{array}
 \right)}$\ \
\rr{N8} ${\scriptsize\left(
 \begin{array}{ccc}
 0&   0&    0\\
 1&  0&   0\\
 0 & 0 & 1
 \end{array}
 \right)}$\ \

 \smallskip

 \rr{N9} ${\scriptsize\left(
 \begin{array}{ccc}
 0&   0&   0\\
 1&  0&   1\\
 0 & 0& 1
 \end{array}
 \right)}$\ \
 \rr{N10} ${\scriptsize\left(
 \begin{array}{ccc}
 0&   0&   1\\
 0&  0&   0\\
 0 & 1& 0
 \end{array}
 \right)}$\ \
 \rr{N11} ${\scriptsize\left(
 \begin{array}{ccc}
 0&   0&   1\\
 1&  0&   0\\
 0 & 0& 0
 \end{array}
 \right)}$\ \
\rr{N12} ${\scriptsize\left(
 \begin{array}{ccc}
 0&   0&   1\\
 0&  0&   0\\
 1 & 0& 0
 \end{array}
 \right)}$\ \

\smallskip

 \rr{N13} ${\scriptsize\left(
 \begin{array}{ccc}
 0&   0&   1\\
 1&  0&   0\\
 1 & 0& 0
 \end{array}
 \right)}$\ \ \
{\rm (O1)} ${\scriptsize\left(
 \begin{array}{ccc}
 1&0&0\\
 0&0& 0\\
 0&0& 0
 \end{array}
 \right)}$\ \
 {\rm (O2)} ${\scriptsize\left(
 \begin{array}{ccc}
 0 &0 &0\\
 1&0& 0\\
 0&0& 0
 \end{array}
 \right)}$\ \
{\rm (O3)} ${\scriptsize\left(
 \begin{array}{ccc}
 0&0&0\\
 1& 0& 0\\
 1&0& 0
 \end{array}
 \right)}$
\end{thm}

\begin{thm}\label{1=2=3 for p=2}
 Let $G$ be a finite $2$-group with $d(G)=3$, $\Phi(G)\le Z(G)$ and $G'\cong C_2^3$. If the type of $G/G'$ is $(2^{m},2^{m},2^{m})$, where $m\ge 2$, then we can choose suitable generators of $G$, such that the characteristic matrix of $G$ is one of the
following matrices, where different types of matrix give non-isomorphic groups.

\medskip

{\rm (P1)} ${\scriptsize\left(
 \begin{array}{ccc}
 1& 0&0\\
 0& 1&0\\
 0& 0&1
 \end{array}
 \right)}$\ \
  {\rm (P2)} ${\scriptsize\left(
 \begin{array}{ccc}
 1& 0&0\\
 0& 1&0\\
 0& 1&1
 \end{array}
 \right)}$\ \
  {\rm (P3)} ${\scriptsize\left(
 \begin{array}{ccc}
 1& 1&1\\
 0& 1&0\\
 0& 0&1
 \end{array}
 \right)}$\ \
  {\rm (P4)} ${\scriptsize\left(
 \begin{array}{ccc}
 1& 0&1\\
 0& 0&1\\
 0& 1&0
 \end{array}
 \right)}$

 \smallskip

 {\rm (Q1)} ${\scriptsize\left(
 \begin{array}{ccc}
 1& 0&0\\
 0& 1&0\\
 0& 0&0
 \end{array}
 \right)}$\ \
  {\rm (Q2)} ${\scriptsize\left(
 \begin{array}{ccc}
 0& 1&0\\
 1& 0&0\\
 0& 0&0
 \end{array}
 \right)}$\ \
  {\rm (Q3)} ${\scriptsize\left(
 \begin{array}{ccc}
 1& 0&0\\
 1& 1&0\\
 0& 0&0
 \end{array}
 \right)}$\ \
  {\rm (Q4)} ${\scriptsize\left(
 \begin{array}{ccc}
 0& 0&1\\
 0& 1&0\\
 0& 0&0
 \end{array}
 \right)}$\ \
 {\rm (Q5)} ${\scriptsize\left(
 \begin{array}{ccc}
 0& 1&0\\
 0& 0&1\\
 0& 0&0
 \end{array}
 \right)}$

 \smallskip

 {\rm (R1)} ${\scriptsize\left(
 \begin{array}{ccc}
 0& 1&0\\
 0& 0&0\\
 0& 0&0
 \end{array}
 \right)}$
 {\rm (R2)} ${\scriptsize\left(
 \begin{array}{ccc}
 1& 0&0\\
 0& 0&0\\
 0& 0&0
 \end{array}
 \right)}$
 {\rm (R3)} ${\scriptsize\left(
 \begin{array}{ccc}
 0& 0&0\\
 0& 0&0\\
 0& 0&0
 \end{array}
 \right)}$
\end{thm}

\begin{thm}\label{N11-property} Let $G$ be a finite $2$-group determined by the characteristic matrix {\rm (N11)} in Theorem $\ref{1>2=3 for m_2=1 and p=2}.$
Then $I_{\min}=2$ and $I_{\max}=m_1+2$.
\end{thm}
\demo Since $\left(
 \begin{array}{cc}
 {w}_{22}& {w}_{23}\\
 {w}_{32} &{w}_{33}
 \end{array}
 \right)=0$,
 by Theorem \ref{maximal-index} (1), $I_{\max}=m_1+2$.
 Since $\rank(w(G))=2$, by Theorem \ref{minimal-index} (2), $1\le I_{\min}\le 2$. We will show that $I_{\min}\neq 1$. Otherwise, by Theorem \ref{minimal-index} (1) and Theorem \ref{isomorphic-3}, there exist invertible matrices $X={\scriptsize\left(
 \begin{array}{ccc}
 1& x_{12}&x_{13}\\
 0 & x_{22}&x_{23}\\
 0 & x_{32}&x_{33}
 \end{array}
 \right)}$ and $X_2={\scriptsize\left(
 \begin{array}{ccc}
 1& 0&0\\
 x_{21}& x_{22}&x_{23}\\
 x_{31}& x_{32}&x_{33}
 \end{array}
 \right)}$, and $w(\bar{G})=X_2w(G)X^t+{\scriptsize\left(
 \begin{array}{ccc}
 0& 0&0\\
 x_{22}x_{23}& 0&0\\
 x_{32}x_{33}&0&0
 \end{array}
 \right)}$ such that $\left(
 \begin{array}{cc}
 \bar{w}_{11}& \bar{w}_{12}\\
 \bar{w}_{21} &\bar{w}_{22}
 \end{array}
 \right)$ is invertible. By calculation, we have
 $W=\left(
 \begin{array}{cc}
 \bar{w}_{11}& \bar{w}_{12}\\
 \bar{w}_{21} &\bar{w}_{22}
 \end{array}
 \right)=\left(
 \begin{array}{cc}
 x_{13}& x_{23}\\
 x_{22}+x_{13}x_{21}+x_{22}x_{23} &x_{21}x_{23}
 \end{array}
 \right)$. Hence $\det(W)=0$, a contradiction. So $I_{\min}=2$.
\qed

\begin{thm}\label{1>2=3==1-property} Let $G$ be a finite $p$-group determined by a characteristic matrix in Theorem $\ref{1>2=3 for m_2=1 and p=2}.$

\rr{1} If $I_{\min}=1$, then the characteristic matrix is one of the following: {\rm (M2)--(M7), (N5)--(N6), (N10), (N12)--(N13)};

\rr{2} If $I_{\min}=2$, then the characteristic matrix is one of the following: {\rm (M1), (N1)--(N4), (N7)--(N9), (N11), (O1)--(O3)};

\rr{3} If $I_{\max}={m_1+2}$, then the characteristic matrix is one of the following: {\rm (N11)--(N13), (O1)--(O3)};

\rr{4} If $I_{\max}={m_1+1}$, then the characteristic matrix is one of the following: {\rm (M4)--(M7), (N1) for $m_1=2$, (N2) for $m_1=2$, (N4)--(N10)};

\rr{5} If $I_{\max}={m_1}$, then the characteristic matrix is one of the following: {\rm (N1) for $m_1>2$, (N2) for $m_1>2$, (N3), (M1)--(M3)}.

\end{thm}

\begin{thm}\label{1=2=3 for p=2-property} Let $G$ be a finite $p$-group determined by a characteristic matrix in Theorem $\ref{1=2=3 for p=2}.$

\rr{1} If $I_{\min}=m$, then the characteristic matrix is one of the following: {\rm (P1)--(P4), (Q1)--(Q5)};

\rr{2} If $I_{\min}=m+1$, then the characteristic matrix is one of the following: {\rm (R1)--(R2)};

\rr{3} If $I_{\min}={m+2}$, then the characteristic matrix is {\rm (R3)};

\rr{4} If $I_{\max}={m+2}$, then the characteristic matrix is one of the following: {\rm (Q1)--(Q2), (Q5), (R1)--(R3)};

\rr{5} If $I_{\max}={m+1}$, then the characteristic matrix is one of the following: {\rm (P1)--(P4), (Q3)--(Q4)}.

\end{thm}

\begin{thm}\label{1=2=3==1} Let $G$ be a finite $2$-group with $d(G)=3$, $\Phi(G)\le Z(G)$ and $G'\cong C_2^3$.
If $|G|=2^6$, then $G$ is one
of the following non-isomorphic groups:

\rr{S1} $\langle a, b, c,d,e,f \di a^{2}=b^{2}=c^{2}=d^2=e^2=f^2=1,
[b,c]=d,[c,a]=e,[a,b]=f,[d,a]=[d,b]=[d,c]=[e,a]=[e,b]=[e,c]=[f,a]=[f,b]=[f,c]=1\rangle$;

\rr{S2} $\langle a, b,c,d,e \di a^{2}=b^{2}=c^{4}=d^2=e^2=1,
[b,c]=d,[c,a]=e,[a,b]=c^2,[d,a]=[d,b]=[d,c]=[e,a]=[e,b]=[e,c]=1\rangle$;

\rr{S3} $\langle a, b,c,d \di a^{4}=b^{4}=c^{2}=d^2=1,
[b,c]=d,[c,a]=a^2b^2,[a,b]=b^2,[d,a]=[d,b]=[d,c]=1\rangle$;

\rr{S4} $\langle a, b,c,d \di a^{4}=b^{4}=c^{2}=d^2=1,
[b,c]=d,[c,a]=a^2b^2,[a,b]=a^2,[d,a]=[d,b]=[d,c]=1\rangle$;

\rr{S5} $\langle a, b,c,d,e \di a^{4}=b^{4}=c^{2}=d^2=e^2=1,
[b,c]=d,[c,a]=e,[a,b]=a^2=b^2,[d,a]=[d,b]=[d,c]=[e,a]=[e,b]=[e,c]=1\rangle$;

\rr{S6} $\langle a, b,c,d,e \di a^{4}=b^{4}=c^{4}=d^2=e^2=1,
[b,c]=d,[c,a]=e,[a,b]=a^2=b^2=c^2,[d,a]=[d,b]=[d,c]=[e,a]=[e,b]=[e,c]=1\rangle$;

\rr{S7} $\langle a, b,c,d \di a^{4}=b^{2}=c^{4}=d^2=1,
[b,c]=d,[c,a]=a^2,[a,b]=c^2,[d,a]=[d,b]=[d,c]=1\rangle$;

\rr{S8} $\langle a,b,c,d \di a^{4}=b^{4}=c^{4}=d^2=1,
[b,c]=d,[c,a]=a^2,[a,b]=b^2=c^2,[d,a]=[d,b]=[d,c]=1\rangle$;

\rr{S9} $\langle a,b,c,d \di a^{4}=b^{4}=c^{4}=d^2=1,
[b,c]=d,[c,a]=a^2b^2,[a,b]=a^2=c^2,[d,a]=[d,b]=[d,c]=1\rangle$;

\rr{S10} $\langle a,b,c \di a^{4}=b^{4}=c^{4}=1,
[b,c]=a^2b^2,[c,a]=b^2c^2,[a,b]=c^2,[c^2,a]=[c^2,b]=1\rangle$;
\end{thm}

\begin{thm}\label{1=2=3==1-property} Let $G$ be a finite $2$-group in Theorem $\ref{1=2=3==1}.$ Then $I_{\min}=1$.

\rr{1} If $I_{\max}={3}$, then $G$ is one of the following: {\rm (S1)--(S2), (S5)--(S6)};

\rr{2} If $I_{\max}=2$, then $G$ is one of the following: {\rm (S3)--(S4), (S7)--(S9)};

\rr{3} If $I_{\max}=1$, then $G$ is {\rm (S10)}.
\end{thm}

\section{The application of the classification}

Before applying the above results, we want to ensure that the classification is accurate. A useful check is to compare determination of the groups whose order is as small as possible. Determination may be the database of Magma or group list given in a paper. For example, if $m_1=m_2=m_3$ and $p$ is odd, then the smallest order of $G$ is $p^6$. In this case, we can compare determinations against both those in Magma and those of James \cite{J}. If $p=2$, then the smallest order of $G$ is $2^9$ for $m_1>m_2>m_3$, $2^{10}$ for $m_1>m_2=m_3>1$ and $2^8$ for $m_1=m_2>m_3$ respectively. In fact, for $m_1>m_2=m_3>1$ and $m_1=m_2>m_3$, we only need to check one of them. Fortunately, Magma provide group lists for 2-groups of order less than $2^{10}$. Hence we can ensure the accuracy for $p=2$.

The above classification can be used easily in classifying $\mathcal{A}_3$-groups and $p$-groups which contain an $\mathcal{A}$-subgroup of index $p$
by investigating the value of $I_{\max}$ and $I_{\min}$ respectively. In this section, we give another application. That is, we pick up all metahamiltonian groups from our results.

\begin{thm}\label{metahamiltonian}
  Suppose that $G$ is a finite $p$-group with $d(G)=3$, $\Phi(G)\le Z(G)$ and $G'\cong C_p^3$. Then $G$ is a metahamiltonian group if and only if $G$ is one of the
following non-isomorphic groups:

\rr{1} $\langle a_1, a_2, a_3 \di a_1^{p^{m_1+1}}=a_2^{p^{m_2+1}}=a_3^{p^{m_3+1}}=1,
[a_2,a_3]=a_1^{p^{m_1}},[a_1,a_3]=a_2^{\eta p^{m_2}},[a_1,a_2]=a_3^{ p^{m_3}}\rangle$,
where $p$ is odd, $m_1=m_2+1=m_3+1$ and $\eta$ is a fixed square
non-residue modulo $p$;

\rr{2} $\langle a_1, a_2, a_3 \di a_1^{p^{m_1+1}}=a_2^{p^{m_2+1}}=a_3^{p^{m_3+1}}=1,
[a_2,a_3]=a_1^{p^{m_1}},[a_1,a_3]=a_2^{l p^{m_2}}a_3^{- p^{m_2}},[a_1,a_2]=a_3^{ p^{m_3}}\rangle$,
where $p$ is odd, $m_1=m_2+1=m_3+1$ and $1+4l\not\in (F_p)^2$;

\rr{3} $\langle a_1, a_2, a_3 \di a_1^{2^{m_1+1}}=a_2^{2^{m_2+1}}=a_3^{2^{m_3+1}}=1,
[a_2,a_3]=a_1^{2^{m_1}},[a_3,a_1]=a_2^{2^{m_2}},[a_1,a_2]=a_2^{2^{m_2}}a_3^{2^{m_3}}\rangle$,
where $m_1=m_2+1=m_3+1$;

\rr{4} $\langle a_1, a_2, a_3 \di a_1^{p^{m_1+1}}=a_2^{p^{m_2+1}}=a_3^{p^{m_3+1}}=1,
[a_2,a_3]=a_1^{p^{m_1}},[a_1,a_3]=a_2^{\eta p^{m_2}},[a_1,a_2]=a_3^{ p^{m_3}}\rangle$,
where $p$ is odd, $m_1=m_2=m_3+1$ and $\eta$ is a fixed square
non-residue modulo $p$;

\rr{5} $\langle a_1, a_2, a_3 \di a_1^{p^{m_1+1}}=a_2^{p^{m_2+1}}=a_3^{p^{m_3+1}}=1,
[a_2,a_3]=a_1^{p^{m_1}},[a_1,a_3]=a_1^{ p^{m_1}}a_2^{l p^{m_2}},[a_1,a_2]=a_3^{ p^{m_3}}\rangle$,
where $p$ is odd, $m_1=m_2=m_3+1$ and $1+4l\not\in (F_p)^2$;

\rr{6} $\langle a_1, a_2, a_3 \di a_1^{2^{m_1+1}}=a_2^{2^{m_2+1}}=a_3^{2^{m_3+1}}=1,
[a_2,a_3]=a_1^{2^{m_1}}a_2^{2^{m_2}} ,[a_3,a_1]=a_2^{2^{m_2}},[a_1,a_2]=a_3^{2^{m_3}}\rangle$,
where $m_1=m_2=m_3+1$;

\rr{7} $\langle a,b,c \di a^{4}=b^{4}=c^{4}=1,
[b,c]=a^2b^2,[c,a]=b^2c^2,[a,b]=c^2,[c^2,a]=[c^2,b]=1\rangle$;

\end{thm}
\demo For most cases, it is easy to find an $\mathcal{A}_1$-subgroup which is not normal in $G$. In fact, in \cite{AF}, we proved that $G'$ contains in every $\mathcal{A}_1$-subgroups for metahamilton $p$-groups. Let $G$ be a metahamiltonian $p$-group and the type of $G/G'$ be $(p^{m_1},p^{m_2},p^{m_3})$. Then we have $I_{\min}=m_3$ and $I_{\max}=m_1$. Hence we only need to check the groups determined by the following characteristic matrices: (A1), (A2), (D2)--(D4), (D6), (D7), (G3) for $-\nu\not\in (F_p^*)^2$, (G4) for $-r\not\in (F_p^*)^2$, (G7), (M2), (M3), and the group (S10).

Case 1: $G$ is a group determined by the characteristic matrix (A1).

In this case, $G=\langle a_1, a_2, a_3 \di a_1^{p^{m_1+1}}=a_2^{p^{m_2+1}}=a_3^{p^{m_3+1}}=1,
[a_2,a_3]=a_1^{p^{m_1}},[a_3,a_1]=a_2^{\nu_1 p^{m_2}},[a_1,a_2]=a_3^{\nu_2 p^{m_3}}\rangle$,
where $m_1>m_2>m_3$, $\nu_1,\nu_2=1$ or a fixed square
non-residue modulo odd $p$. Since $\lg a_2,a_3a_1^p\rg$ is neither abelian nor normal in $G$, $G$ is not a metahamiltonian group.

Case 2: $G$ is a group determined by the characteristic matrix (A2).

In this case, $G=\langle a_1, a_2, a_3 \di a_1^{p^{m_1+1}}=a_2^{p^{m_2+1}}=a_3^{p^{m_3+1}}=1,
[a_2,a_3]=a_1^{p^{m_1}},[a_3,a_1]^t=a_3^{ p^{m_3}},[a_1,a_2]=a_2^{ p^{m_2}}\rangle$,
where $m_1>m_2>m_3$ and $t \neq 0$. Since $\lg a_1,a_2\rg$ is neither abelian nor normal in $G$, $G$ is not a metahamiltonian group.

Case 3: $G$ is a group determined by the characteristic matrix (D2).

In this case, $G=\langle a_1, a_2, a_3 \di a_1^{p^{m_1+1}}=a_2^{p^{m_2+1}}=a_3^{p^{m_3+1}}=1,
[a_2,a_3]=a_1^{p^{m_1}},[a_3,a_1]=a_3^{-p^{m_3}},[a_1,a_2]=a_2^{ p^{m_2}}a_3^{\nu p^{m_3}}\rangle$,
where $p$ is odd, $m_1>m_2=m_3$, $\nu=1$ or a fixed square
non-residue modulo $p$. Since $\lg a_1,a_2\rg$ is neither abelian nor normal in $G$, $G$ is not a metahamiltonian group.

Case 4: $G$ is a group determined by the characteristic matrix (D3).

In this case, $G=\langle a_1, a_2, a_3 \di a_1^{p^{m_1+1}}=a_2^{p^{m_2+1}}=a_3^{p^{m_3+1}}=1,
[a_2,a_3]=a_1^{p^{m_1}},[a_3,a_1]=a_2^{ p^{m_2}},[a_1,a_2]^\nu=a_3^{ p^{m_3}}\rangle$,
where $p$ is odd, $m_1>m_2=m_3$, $\nu=1$ or a fixed square
non-residue modulo $p$. If $m_1>m_2+1$, then $\lg a_2,a_3a_1^p\rg$ is neither abelian nor normal in $G$. Hence $G$ is not a metahamiltonian group. In the following, we may assume that $m_1=m_2+1=m_3+1$.

Notice that $G$ is not a metahamiltonian group if and only if there exists $D\in\mathcal{A}_1$ such that $G'\not\le D$. Since $|G:D|\le p^{m_3+1}$, $|G:DG'|=|G:D|/|DG':D|\le p^{m_3}$. On the other hand, since $DG'/G'$ is a subgroup of $G/G'$ and $d(DG'/G')=2$, we have $|G/G':DG'/G'|\ge p^{m_3}$. It follows that $|G:DG'|=p^{m_3}$ and hence we may assume that $D=\lg \bar{a}_1,\bar{a}_2\rg$. By Theorem \ref{isomorphic-1}, there exist invertible matrices $X_2={\scriptsize\left(
 \begin{array}{ccc}
 x_{11}&      0&      0\\
 x_{21}& x_{22}&      x_{23}\\
 x_{31}& x_{32}&   x_{33}
 \end{array}
 \right)}$ and $X={\scriptsize\left(
 \begin{array}{ccc}
x_{11} & x_{12}  & x_{13}\\
0 & x_{22}  &x_{23}\\
0 &  x_{32}  & x_{33}
 \end{array}
 \right)}$, and
$w(\bar{G})=\det(X)^{-1}X_2w(G)X^t$ such that $\rank{\scriptsize\left(
 \begin{array}{ccc}
\bar{w}_{11} & \bar{w}_{12}  & \bar{w}_{13}\\
\bar{w} & \bar{w}_{22}  &\bar{w}_{23}\\
0 &  0  & 1
 \end{array}
 \right)}=2$.
By calculation, we have $\det(X)w(\bar{G})={\scriptsize\left(
 \begin{array}{ccc}
x_{11}^2 & 0  & *\\
* & x_{22}^2+\nu x_{23}^2&*\\
* &  *  & *
 \end{array}
 \right)}$. Hence $G$ is not a metahamiltonian group if and only if there exist $x_{22}$ and $x_{33}$ such that $x_{22}^2+\nu x_{23}^2=0$. That is, $-\nu\in (F_p)^2$. On another words, $G$ is a metahamiltonian group if and only if $-\nu\not\in (F_p)^2$.

Let $\eta$ be a fixed square
non-residue modulo $p$. If $-\nu\not\in (F_p^*)^2$, then there exists $s\in F_p$ such that $-\nu\eta= s^2$. Replacing $a_1,a_2$ and $a_3$ with $a_1^s,a_3^{-1}$ and $a_2^s$ respectively, we get the group of Type (1).

Case 5: $G$ is a group determined by the characteristic matrix (D4).

In this case, $G=\langle a_1, a_2, a_3 \di a_1^{p^{m_1+1}}=a_2^{p^{m_2+1}}=a_3^{p^{m_3+1}}=1,
[a_2,a_3]=a_1^{p^{m_1}},[a_3,a_1]^{r+1}=a_2^{r p^{m_2}}a_3^{-p^{m_3}},[a_1,a_2]^{r+1}=a_2^{p^{m_2}}a_3^{ p^{m_3}}\rangle$,
where $p$ is odd, $m_1>m_2=m_3$, $r=1,2,\dots,p-2$. If $m_1>m_2+1$, then $\lg a_2,a_3a_1^p\rg$ is neither abelian nor normal in $G$. Hence $G$ is not a metahamiltonian group. In the following, we may assume that $m_1=m_2+1=m_3+1$.

Notice that $G$ is not a metahamiltonian group if and only if there exists $D\in\mathcal{A}_1$ such that $G'\not\le D$. Since $|G:D|\le p^{m_3+1}$, $|G:DG'|=|G:D|/|DG':D|\le p^{m_3}$. On the other hand, since $DG'/G'$ is a subgroup of $G/G'$ and $d(DG'/G')=2$, we have $|G/G':DG'/G'|\ge p^{m_3}$. It follows that $|G:DG'|=p^{m_3}$ and hence we may assume that $D=\lg \bar{a}_1,\bar{a}_2\rg$. By Theorem \ref{isomorphic-1}, there exist invertible matrices $X_2={\scriptsize\left(
 \begin{array}{ccc}
 x_{11}&      0&      0\\
 x_{21}& x_{22}&      x_{23}\\
 x_{31}& x_{32}&   x_{33}
 \end{array}
 \right)}$ and $X={\scriptsize\left(
 \begin{array}{ccc}
x_{11} & x_{12}  & x_{13}\\
0 & x_{22}  &x_{23}\\
0 &  x_{32}  & x_{33}
 \end{array}
 \right)}$, and
$w(\bar{G})=\det(X)^{-1}X_2w(G)X^t$ such that $\rank{\scriptsize\left(
 \begin{array}{ccc}
\bar{w}_{11} & \bar{w}_{12}  & \bar{w}_{13}\\
\bar{w} & \bar{w}_{22}  &\bar{w}_{23}\\
0 &  0  & 1
 \end{array}
 \right)}=2$.
By calculation, we have $\det(X)w(\bar{G})={\scriptsize\left(
 \begin{array}{ccc}
x_{11}^2 & 0  & 0\\
* & x_{22}^2+r x_{23}^2&*\\
* &  *  & *
 \end{array}
 \right)}$. Hence $G$ is not a metahamiltonian group if and only if there exist $x_{22}$ and $x_{33}$ such that $x_{22}^2+r x_{23}^2=0$. That is, $-r\in (F_p)^2$. On another words, $G$ is a metahamiltonian group if and only if $-r\not\in (F_p)^2$. Now we assume that $G$ is a metahamiltonian group.

If $(r+1)r^{-1}\in (F_p^*)^2$, then, letting $(r+1)r^{-1}=s^2$, and replacing $a_1, a_2$ and $a_3$ with $a_1^{-2^{-1}(r+1)}$, $a_3^{-s}$ and $a_2^{-2^{-1}rs}a_3^{2^{-1}s}$ respectively, we get the group of Type (2).

If $r\in (F_p^*)^2$, then, letting $r=s^2$, and replacing $a_1, a_2$ and $a_3$ with $a_1^{-2^{-1}(r+1)}$, $a_2^{-s}a_3^{s^{-1}}$, and $a_2^{s}a_3^{2^{-1}s^{-1}(r-1)}$ respectively, we get the group of Type (2).

If $r+1\in (F_p^*)^2$, then, letting $r+1=s^2$, and replacing $a_1, a_2$ and $a_3$ with $a_1^{-2^{-1}(r+1)}$, $a_2^{s-2s^{-1}}a_3^{-2s^{-1}}$, and $a_2^{2s^{-1}-2^{-1}3s}a_3^{2s^{-1}-2^{-1}s}$ respectively, we get the group of Type (2).

Since there must exist a square among $r,r+1$ and $(r+1)r^{-1}$, we get the group of Type (2). For all cases, we have $l=-4^{-1}(r+1)$. Since $-r\not\in (F_p)^2$, we have $1+4l\not\in (F_p)^2$.

Case 6: $G$ is a group determined by the characteristic matrix (D6).

In this case, $G=\langle a_1, a_2, a_3 \di a_1^{2^{m_1+1}}=a_2^{2^{m_2+1}}=a_3^{2^{m_3+1}}=1,
[a_2,a_3]=a_1^{2^{m_1}},[a_3,a_1]=a_2^{2^{m_2}},[a_1,a_2]=a_2^{2^{m_2}}a_3^{2^{m_3}}\rangle$,
where $m_1>m_2=m_3>1$. Since $\lg a_1,a_2a_3\rg$ is neither abelian nor normal in $G$, $G$ is not a metahamiltonian group.

Case 7: $G$ is a group determined by the characteristic matrix (D7).

In this case, $G=\langle a_1, a_2, a_3 \di a_1^{2^{m_1+1}}=a_2^{2^{m_2+1}}=a_3^{2^{m_3+1}}=1,
[a_2,a_3]=a_1^{2^{m_1}},[a_3,a_1]=a_2^{2^{m_2}},[a_1,a_2]=a_2^{2^{m_2}}a_3^{2^{m_3}}\rangle$,
where $m_1>m_2=m_3>1$. If $m_1>m_3+1$, then $\lg a_2,a_3a_1^2\rg$ is neither abelian nor normal in $G$. Hence $G$ is not a metahamiltonian group. If $m_1=m_2+1$, then, in the following, we will prove that $G$ is a metahamiltonian group.

 If $m_1=m_2+1$, then we claim that $G'$ contains in every $\mathcal{A}_1$-subgroups of $G$. Otherwise, there exists $D\in\mathcal{A}_1$ such that $G'\not\le D$. Since $|G:D|\le 2^{m_3+1}$, $|G:DG'|=|G:D|/|DG':D|\le 2^{m_3}$. On the other hand, since $DG'/G'$ is a subgroup of $G/G'$ and $d(DG'/G')=2$, we have $|G/G':DG'/G'|\ge 2^{m_3}$. It follows that $|G:DG'|=2^{m_3}$ and hence we may assume that $D=\lg a_1a_2^ib, a_3c\rg$ or $D=\lg a_1a_3^jb,a_2a_3^kc\rg$ where $b,c\in \Phi(G)$. It is easy to see, by calculation, that $G'\le D$, a contradiction. Hence $G'$ contains in every $\mathcal{A}_1$-subgroups of $G$, and $G$ is a metahamiltonian group. In this case, we get the group of Type (3) where $m_3>1$.

Case 8: $G$ is a group determined by the characteristic matrix (G3) for $-\nu\not\in (F_p^*)^2$.

In this case, $G=\langle a_1, a_2, a_3 \di a_1^{p^{m_1+1}}=a_2^{p^{m_2+1}}=a_3^{p^{m_3+1}}=1,
[a_2,a_3]=a_1^{p^{m_1}},[a_3,a_1]^\nu=a_2^{ p^{m_2}},[a_1,a_2]=a_3^{ p^{m_3}}\rangle$,
where $p$ is odd, $m_1=m_2>m_3$, $\nu=1$ or a fixed square
non-residue modulo $p$, and $-\nu\not\in (F_p^*)^2$. If $m_2>m_3+1$, then $\lg a_2,a_3a_1^p\rg$ is neither abelian nor normal in $G$. Hence $G$ is not a metahamiltonian group. If $m_2=m_3+1$, then, in the following, we will prove that $G$ is a metahamiltonian group.

 If $m_2=m_3+1$, then we claim that $G'$ contains in every $\mathcal{A}_1$-subgroups of $G$. Otherwise, there exists $D\in\mathcal{A}_1$ such that $G'\not\le D$. Since $|G:D|\le p^{m_3+1}$, $|G:DG'|=|G:D|/|DG':D|\le p^{m_3}$. On the other hand, since $DG'/G'$ is a subgroup of $G/G'$ and $d(DG'/G')=2$, we have $|G/G':DG'/G'|\ge p^{m_3}$. It follows that $|G:DG'|=p^{m_3}$ and hence we may assume that $D=\lg a_1a_3^jb,a_2a_3^kc\rg$ where $b,c\in \Phi(G)$. It is easy to see, by calculation, that $G'\le D$, a contradiction. Hence $G'$ contains in every $\mathcal{A}_1$-subgroups of $G$, and $G$ is a metahamiltonian group.

Let $\eta$ be a fixed square
non-residue modulo $p$. Since $-\nu\not\in (F_p^*)^2$, there exists $s\in F_p$ such that $-\nu\eta= s^2$. Replacing $a_1$ and $a_3$ with $a_1^s$ and $a_3^s$ respectively, we get the group of Type (4).

Case 9: $G$ is a group determined by the characteristic matrix (G4) for $-r\not\in (F_p^*)^2$.

In this case, $G=\langle a_1, a_2, a_3 \di a_1^{p^{m_1+1}}=a_2^{p^{m_2+1}}=a_3^{p^{m_3+1}}=1,
[a_2,a_3]^{1+r}=a_1^{rp^{m_1}}a_2^{p^{m_2}},[a_3,a_1]^{r+1}$ $=a_1^{-p^{m_1}}a_2^{p^{m_2}},[a_1,a_2]=a_3^{ p^{m_3}}\rangle$,
where $p$ is odd, $m_1=m_2>m_3$, $r=1,2,\dots,p-2$, and $-r\not\in (F_p^*)^2$. If $m_2>m_3+1$, then $\lg a_2,a_3a_1^p\rg$ is neither abelian nor normal in $G$. Hence $G$ is not a metahamiltonian group. If $m_2=m_3+1$, then, in the following, we will prove that $G$ is a metahamiltonian group.

 If $m_2=m_3+1$, then we claim that $G'$ contains in every $\mathcal{A}_1$-subgroups of $G$. Otherwise, there exists $D\in\mathcal{A}_1$ such that $G'\not\le D$. Since $|G:D|\le p^{m_3+1}$, $|G:DG'|=|G:D|/|DG':D|\le p^{m_3}$. On the other hand, since $DG'/G'$ is a subgroup of $G/G'$ and $d(DG'/G')=2$, we have $|G/G':DG'/G'|\ge p^{m_3}$. It follows that $|G:DG'|=p^{m_3}$ and hence we may assume that $D=\lg a_1a_3^jb,a_2a_3^kc\rg$ where $b,c\in \Phi(G)$. It is easy to see, by calculation, that $G'\le D$, a contradiction. Hence $G'$ contains in every $\mathcal{A}_1$-subgroups of $G$, and $G$ is a metahamiltonian group.

If $(r+1)r^{-1}\in (F_p^*)^2$, then, letting $(r+1)r^{-1}=s^2$, and replacing $a_1, a_2$ and $a_3$ with $a_1^{2^{-1}rs}a_2^{2^{-1}s}$, $a_2^{s}$ and $a_3^{2^{-1}(1+r)}$ respectively, we get the group of Type (5).

If $r\in (F_p^*)^2$, then, letting $r=s^2$, and replacing $a_1, a_2$ and $a_3$ with $a_2^{-2^{-1}s^{-1}(r+1)}$, $a_1^{s}a_2^{-s^{-1}}$, and $a_3^{2^{-1}(r+1)}$ respectively, we get the group of Type (5).

If $r+1\in (F_p^*)^2$, then, letting $r+1=s^2$, and replacing $a_1, a_2$ and $a_3$ with $(a_1a_2)^{2^{-1}s}$, $(a_1a_2)^{2s^{-1}}a_1^{-s}$, and $a_3^{2^{-1}(r+1)}$ respectively, we get the group of Type (5).

Since there must exist a square in $r,r+1$ and $(r+1)r^{-1}$, we get the group of Type (5). For all cases, we have $l=-4^{-1}(r+1)$. Since $-r\not\in (F_p)^2$, we have $1+4l\not\in (F_p)^2$.

Case 10: $G$ is a group determined by the characteristic matrix (G7).

In this case, $G=\langle a_1, a_2, a_3 \di a_1^{2^{m_1+1}}=a_2^{2^{m_2+1}}=a_3^{2^{m_3+1}}=1,
[a_2,a_3]=a_1^{2^{m_1}}a_2^{2^{m_2}} ,[a_3,a_1]=a_2^{2^{m_2}},[a_1,a_2]=a_3^{2^{m_3}}\rangle$,
where $m_1=m_2>m_3$. If $m_2>m_3+1$, then $\lg a_2,a_3a_1^2\rg$ is neither abelian nor normal in $G$. Hence $G$ is not a metahamiltonian group. If $m_2=m_3+1$, then, in the following, we will prove that $G$ is a metahamiltonian group.

 If $m_2=m_3+1$, then we claim that $G'$ contains in every $\mathcal{A}_1$-subgroups of $G$. Otherwise, there exists $D\in\mathcal{A}_1$ such that $G'\not\le D$. Since $|G:D|\le 2^{m_3+1}$, $|G:DG'|=|G:D|/|DG':D|\le 2^{m_3}$. On the other hand, since $DG'/G'$ is a subgroup of $G/G'$ and $d(DG'/G')=2$, we have $|G/G':DG'/G'|\ge 2^{m_3}$. It follows that $|G:DG'|=2^{m_3}$ and hence we may assume that $D=\lg a_1a_3^jb,a_2a_3^kc\rg$ where $b,c\in \Phi(G)$. It is easy to see, by calculation, that $G'\le D$, a contradiction. Hence $G'$ contains in every $\mathcal{A}_1$-subgroups of $G$, and $G$ is a metahamiltonian group. In this case, we get the group of Type (6).

Case 11: $G$ is a group determined by the characteristic matrix (M2).

In this case, $G=\langle a_1, a_2, a_3 \di a_1^{2^{m_1+1}}=a_2^{4}=a_3^{4}=1,
[a_2,a_3]=a_1^{2^{m_1}} ,[a_3,a_1]=a_2^{2},[a_1,a_2]=a_3^{2}\rangle$,
where $m_1>1$. Since $\lg a_1,a_2a_3\rg$ is neither abelian nor normal in $G$, $G$ is not a metahamiltonian group.

Case 12: $G$ is a group determined by the characteristic matrix (M3).

In this case, $G=\langle a_1, a_2, a_3 \di a_1^{2^{m_1+1}}=a_2^{4}=a_3^{4}=1,
[a_2,a_3]=a_1^{2^{m_1}} ,[a_3,a_1]=a_2^{2},[a_1,a_2]=a_2^2a_3^{2}\rangle$,
where $m_1>1$. If $m_1>2$, then $\lg a_2,a_3a_1^2\rg$ is neither abelian nor normal in $G$. Hence $G$ is not a metahamiltonian group. If $m_1=2$, then, in the following, we will prove that $G$ is a metahamiltonian group.

 If $m_1=2$, then we claim that $G'$ contains in every $\mathcal{A}_1$-subgroups of $G$. Otherwise, there exists $D\in\mathcal{A}_1$ such that $G'\not\le D$. Since $|G:D|\le 2^2$, $|G:DG'|=|G:D|/|DG':D|\le 2$. It follows that $|G:DG'|=2$ and hence we may assume that $D=\lg a_1a_2^ib,a_3c\rg$ or $D=\lg a_1a_3^jb,a_2a_3^kc\rg$ where $b,c\in \Phi(G)$. It is easy to see, by calculation, that $G'\le D$, a contradiction. Hence $G'$ contains in every $\mathcal{A}_1$-subgroups of $G$, and $G$ is a metahamiltonian group. In this case, we get the group of Type (3) where $m_3=1$.

Case 13: $G$ is the group (S10).

Notice that $I_{\max}=1$ for the group (S10). It means all $\mathcal{A}_1$-subgroups are maximal in $G$. Hence $G$ is a metahamiltonian group. In fact, $G$ is the smallest Suzuki 2-group and an $\mathcal{A}_2$-group. In this case, we get the group of Type (7).\qed

\end{document}